\documentclass[a4paper,11pt]{amsart}
\usepackage{amsfonts}
\usepackage{amssymb}
\usepackage[utf8]{inputenc}
\usepackage{amsmath}
\usepackage{pdflscape}
\usepackage{easyReview}
\usepackage{graphicx}
\usepackage[colorlinks=true, linkcolor=red, citecolor=blue]{hyperref}

\usepackage[]{epsfig}
\usepackage[]{pstricks}
\usepackage{tikz}
\usepackage[pagewise]{lineno}

\newtheorem{theorem}{Theorem}[section]
\newtheorem*{theorem*}{Theorem}
\newtheorem{proposition}[theorem]{Proposition}
\newtheorem{lemma}[theorem]{Lemma}

\theoremstyle{remark}

\newtheorem{definition}{Definition}
\usepackage{mathrsfs}

\numberwithin{equation}{section}
\makeatletter
\@namedef{subjclassname@2010}{\textup{2020} Mathematics Subject Classification}
\makeatother

\begin{document}
	
	\pagenumbering{arabic}	
	\title[Boundary observation for the KP-II]{
Critical lengths for the linear Kadomtsev-Petviashvili II equation
}
\author[Capistrano-Filho]{Roberto de A. Capistrano--Filho$^{\star}$}
\address{Departamento de Matem\'atica,  Universidade Federal de Pernambuco (UFPE), 50740-545, Recife (PE), Brazil.}
			\email{roberto.capistranofilho@ufpe.br}

\author[Gallego]{Fernando A. Gallego}
\address{Departamento de Matematica, Universidad Nacional de Colombia (UNAL), 170003, Cra 27 No. 64-60, Manizales, Caldas, Colombia}
\email{fagallegor@unal.edu.co}

\author[Muñoz]{Juan R. Muñoz}
\address{Faculty of Electrical Engineering and Applied Computing, University of Dubrovnik,\\
Dubrovnik, 20000, Croatia}
\address{Departamento de Ingeniería Matemática, Facultad de Ciencias Físicas y Matemáticas, Universidad de Chile, Santiago, 8370456, Chile}
			\email{juan.munoz@dim.uchile.cl}

	\thanks{$^\star$Corresponding author: roberto.capistranofilho@ufpe.br}
	\date{\today}
	
	\begin{abstract}
The critical length phenomenon of the Korteweg-de Vries equation is well known; however, in higher dimensions, it is unknown. This work explores this property in the context of the Kadomtsev-Petviashvili equation, a two-dimensional generalization of the Korteweg-de Vries equation. Specifically, we demonstrate observability inequalities for this equation, which allow us to deduce the exact boundary controllability and boundary exponential stabilization of the linear system, provided that the spatial domain length avoids certain specific values, a direct consequence of the Paley-Wiener theorem. To the best of our knowledge, our work introduces new results by identifying a set of critical lengths for the two-dimensional Kadomtsev-Petviashvili equation.
	\end{abstract}

	\maketitle
	
\section{Introduction}
The renowned Korteweg-de Vries (KdV) equation
$$
u_t +  uu_x + u_{xxx} = 0
$$
was first derived by Boussinesq~\cite{Boussinesq1895} and Korteweg-de Vries~\cite{KdV1895} to model the unidirectional propagation of small-amplitude, long-wavelength water waves in shallow canals. This equation encapsulates the interplay between nonlinear and dispersive effects, which are crucial for the formation and stability of solitary waves, often referred to as solitons. The KdV equation is essential in the study of nonlinear dispersive wave systems, offering profound insights into wave dynamics in one spatial dimension.

In scenarios where wave propagation extends to two spatial dimensions, phenomena characterized by weak transverse effects and weak nonlinearity are introduced by Kadomtsev and Petviashvili~\cite{KP70}, giving the Kadomtsev-Petviashvili (KP) equation:
\begin{equation}\label{eq:i.KP}
    \left(u_t + uu_x + a u_{xxx}\right)_x + b u_{yy} = 0,
\end{equation}
where $u = u(x,y,t)$ is the wave profile, and $a$ and $b$ are real constants. The KP equation generalizes the KdV equation to account for the transverse spatial variable $y$, incorporating weak transverse perturbations into the dynamics. It arises naturally in systems where the longitudinal nonlinear and dispersive effects dominate, while transverse variations introduce secondary corrections.

Notice that, through suitable scaling transformations, the coefficients $a$ and $b$ can be normalized. Specifically, by applying the transformations $x \mapsto a^{1/4}x$ and $y \mapsto |b|^{1/2}y$, the equation simplifies to a form with $a = 1$ and $b^2 = 1$. Consequently, we consider the rescaled and simplified version of \eqref{eq:i.KP} for further analysis:
$$
    u_t + uu_x + u_{xxx} + b \partial_x^{-1}u_{yy} = 0,\quad b = \pm 1.
$$
The parameter $b$ determines the type of KP equation: when $b = -1$, the equation is referred to as KP-I, which is relevant for situations where surface tension effects are dominant, leading to a focusing regime. Conversely, when $b = 1$, it is called KP-II, representing a defocusing regime often associated with gravity-dominated systems. These two cases exhibit distinct mathematical and physical wave behaviors.

\subsection{Problem setting}
Numerous aspects of the KP equation have been thoroughly examined, such as its well-posedness, stability of solitary waves, its integrability, etc (see~\cite{Hadac, Iorio98, Isaza01, Levandosky00, Levandosky01, Panthee05, Takaoka00} and therein). Nonetheless, \cite{Saut2022} gives a complete revision of relevant results for the KP and related equations. However, the controllability problem needs to be studied better.

The controllability problem consists of determining whether a suitable control can steer the solution of a system from a prescribed initial state to a desired final state within a given time. Another fundamental objective of control theory is stabilization, in which the control is designed to modify the long-time behavior of the solutions and, in particular, to force them to converge toward a prescribed equilibrium. 

For the KP-II equation, the internal controllability on a periodic domain was studied in~\cite{Rivas2020}. Regarding stabilization, the asymptotic behavior of solutions for the KP-II and their fifth-order generalization, called the Kawahara-Kadomtsev-Petviashvili-II equation, or shorter, K-KP-II equations with internal feedback mechanisms on bounded domains was examined in~\cite{Panthee2011, Moura2022}, respectively. Additionally, the exponential stabilization of the K-KP-II equation with internal damping and time-delayed feedback was investigated in~\cite{Munoz2023}.

Two important facts need to be spotlighted: 
\begin{enumerate} 
\item Due to the regularity of the solutions for the KP-II equation and the dimension of the domain, the Sobolev embeddings employed in the one-dimensional cases cannot be used to address the well-posedness for the nonlinear case; precisely, as $H^1(\Omega)$ is not embedded in $L^\infty(\Omega)$, the fixed point argument does not hold. Therefore, we deal with the linear version of the KP-II.
\item The term $\partial_x u$, called drift, could modify the behavior of the operator in a bounded rectangular domain. Switching to traveling coordinates is not applicable in the bounded domain problem without incurring a real cost. By changing variables to $u(x, t) = u(x + t, t)$, we can eliminate the problematic term in the evolution equation. However, the boundary condition is now applied at a varying spatial point, and the problem is defined in the unusual domain $\{(x, t) \colon t \geq 0, x + t \geq 0, L\geq x\}$ instead of a traditional rectangular domain.
\end{enumerate}

In this work, we will consider the linear KP-II within a rectangular domain $\Omega := (0,L)\times(0,L)$, $L>0$
\begin{equation}\label{eq:KP}
u_t + u_x + u_{xxx} +\partial_x^{-1}(u_{yy})  = 0, \quad (x,y)\in \Omega ,\ t\in(0,T),
\end{equation} 
with  boundary conditions and initial data as follows
\begin{equation}\label{eq:BC}
\begin{cases}
u(0,y,t) = u(L,y,t) = 0, \quad u_x(L,y,t) = h(y,t), & y\in(0,L) ,\ t\in(0,T), \\
u(x,0,t) = u(x,L,t) = 0, & x\in(0,L) ,\ t\in(0,T),\\
u(x,y,0) = u_0(x,y),& (x,y)\in\Omega.
\end{cases}
\end{equation}
 Here, the operator $\partial_x^{-1}$ is defined by\footnote{It can be shown that the definition of operator $\partial_x^{-1}$ is equivalent to $\partial_x^{-1} u(x,y,t) = -\int_x^L u(s,y,t)\,ds$.}
\begin{align*}
\partial_x^{-1} \varphi(x,y,t)= \psi(x,y,t) \quad \text{such that} \quad  \psi(L,y,t)=0 \quad \text{and} \quad \psi_x(x,y,t)=\varphi(x,y,t).
\end{align*}

The boundary condition in \eqref{eq:BC} introduces a natural boundary input acting through the Neumann trace at $x=L$. This naturally leads to the first fundamental issue in control theory, namely, whether this boundary action is sufficient to steer the system between arbitrary states. More precisely, we investigate the following exact boundary controllability problem.

\vspace{0.2cm}

\noindent\emph{\textbf{Problem A:} Given an initial state $u_0$ and a final state $u_1$ in a certain space, can one find an appropriate boundary control
input  so that the equation~\eqref{eq:KP}-\eqref{eq:BC} admits a solution $u$ which equals $u_0$ at time $t = 0$ and $u_1$ at time $t = T$ ?}

\vspace{0.2cm}

The same boundary input can be used for a different purpose. Instead of prescribing it beforehand to steer the solution from an initial state to a desired final state, one may determine it from the current state of the system. Such a choice is called a boundary feedback law and naturally leads to the stabilization problem.

The natural energy associated with \eqref{eq:KP} is
\begin{equation}\label{eq:S.En}
E(t) = \frac{1}{2} \int_{\Omega} u^2(x,y,t)\,dx\,dy.
\end{equation}
The objective is to design a boundary feedback law that dissipates the total energy of the system, ideally with an explicit decay rate. This raises the following question.

\vspace{0.2cm}

\noindent\emph{\textbf{Problem B:} Is it possible to choose the control $h(y,t)$ as a feedback damping mechanism such that $E(t) \to 0$ as $t\to +\infty?$ If this is the case, can we give the decay rate?}

\vspace{0.2cm}

To the best of our knowledge, the boundary control properties of the KP-II equation posed on a rectangle remain completely open. The purpose of this paper is to address the two problems formulated above. More precisely, we establish suitable boundary observability inequalities for the associated conservative system. These inequalities constitute the key ingredient in proving both the exact boundary controllability and the asymptotic stabilization under a single boundary feedback law.

\subsection{Background}
Controllability and stabilization problems have been extensively investigated for a wide variety of one-dimensional dispersive equations and systems, including the KdV, Kawahara, and Boussinesq equations, posed on different domains and under various control configurations; see, for instance, \cite{Rosier,Glass09,Gallego2018,Capistrano2019} and the references therein. 

In applications, it is often desirable to achieve controllability or stabilization using as few control inputs as possible. In the present setting, we consider the boundary conditions~\eqref{eq:BC}, with only one control acting along the right side of the rectangle $\Omega$. For this control configuration, the interaction between the transport term $\partial_xu$, the dispersive operator, and the imposed boundary conditions plays a fundamental role in the appearance of \emph{critical lengths}. More precisely, the proof of the observability inequality reduces to establishing a suitable unique continuation property for the corresponding adjoint system. The failure of this property is characterized by the existence of nontrivial solutions to an overdetermined spectral problem. The values of the length $L$ for which such nontrivial solutions exist, and consequently exact controllability fails, are called critical lengths.

The critical-length phenomenon was first identified by Rosier~\cite{Rosier} in his pioneering work on the linear KdV equation. More precisely, for the boundary control system
\begin{equation}\label{eq:KdV.control}
\begin{cases}
u_t+u_{xxx}+u_x=0,
    &(x,t)\in(0,L)\times(0,T),\\
u(0,t)=u(L,t)=0,\qquad
u_x(L,t)=h(t),
    &t\in(0,T),
\end{cases}
\end{equation}
the set of critical lengths is explicitly given by
\begin{equation}\label{lionel-rosier}
\mathcal N=
\left\{
\frac{2\pi}{\sqrt{3}}
\sqrt{k^2+kl+l^2};
\quad
k,l\in\mathbb Z^+
\right\}.
\end{equation}

The derivation of this characterization relies on two main ingredients. First, the observability inequality is reduced to a unique continuation property for the adjoint system. Second, the failure of this property is shown to be equivalent to the existence of a nontrivial solution to the overdetermined spectral problem
\begin{equation}\label{eq:KdV.sp}
\begin{cases}
\lambda u+u'+u'''=0,\\
u(0)=u(L)=u'(0)=u'(L)=0.
\end{cases}
\end{equation}
The corresponding spectral characterization is given by the following result.

\begin{lemma}[{\cite[Lemma~3.5]{Rosier}}]
Let $L\in(0,+\infty)$. There exist $\lambda\in\mathbb C$ and a nontrivial solution
$u\in H^3(0,L)$ of~\eqref{eq:KdV.sp}
if and only if $L\in\mathcal N$.
\end{lemma}

Following Rosier's characterization, several authors have investigated the characterization of critical lengths for different dispersive equations and boundary control configurations. One natural extension concerns higher-order dispersive equations. In particular, Araruna, Doronin, and Rosier~\cite{ArarunaKawahara} considered the Kawahara equation posed on a bounded interval with a single boundary control and proved that local exact controllability holds whenever the interval length avoids a countable set of critical lengths. Their analysis extends Rosier's spectral approach from the third-order KdV equation to a fifth-order dispersive model. 

For the KdV equation, the influence of critical lengths on nonlinear boundary controllability has been extensively investigated. Coron and Cr\'epeau~\cite{Crepeau2} proved that, for certain critical lengths at which the corresponding linearized system is not controllable, the nonlinear term restores local exact controllability around the origin. This analysis was subsequently extended by Cerpa~\cite{Cerpa}, who established local exact controllability on critical spatial domains provided that the control time is sufficiently large. Cerpa and Cr\'epeau~\cite{CerpaCrepeau2009} later obtained local controllability for the nonlinear KdV equation on any critical domain, again under an appropriate lower bound on the control time. More recently, Coron, Koenig, and Nguyen~\cite{CoronKoenigNguyen} studied the small-time local controllability problem at critical lengths and proved that, for a class of such lengths, the nonlinear KdV system is not small-time locally controllable. For the KdV equation with a right Dirichlet boundary control, the corresponding critical-length obstruction was also analyzed by Glass and Guerrero~\cite{GG} and, more recently, by Nguyen~\cite{Nguyen1}. These works show that critical lengths affect not only the controllability of the linearized system but also the nonlinear recovery of controllability and the minimal control time. 

Employing an alternative approach, \cite{Xiang} proved the null controllability of a linearized KdV equation via a backstepping approach, providing explicit feedback constructions and constructive stabilization procedures as an alternative to the Hilbert Uniqueness Method. We also mention that critical-length phenomena have also been investigated for multidimensional dispersive equations, including the Zakharov--Kuznetsov equation posed on bounded rectangles and strips~\cite{Doronin}.

Concerning stabilization, several increasingly strong decay properties are now available for the KdV equation. Cerpa and Cr\'epeau~\cite{CerpaCrepeauRapid} established rapid exponential stabilization for the linear KdV equation with homogeneous Dirichlet boundary conditions and a Neumann control at the right endpoint by constructing feedback laws that force the corresponding closed-loop solutions to decay exponentially with an arbitrarily prescribed rate.

For the nonlinear KdV equation with the same right Neumann boundary configuration, Coron and L\"u~\cite{CoronLu} established local rapid stabilization for every noncritical length. Their approach is based on the construction of a suitable integral transformation and yields exponential decay with an arbitrarily large prescribed rate, provided that the initial data is sufficiently small.

Nguyen~\cite{Nguyen} investigated the asymptotic decay and stabilization of nonlinear KdV equations at critical lengths, establishing polynomial decay estimates and asymptotic stability for small initial data, even in regimes where the linearized dynamics preserve the energy on finite-dimensional invariant subspaces. More recently, Nguyen~\cite{NguyenStabilization} obtained both local rapid and finite-time boundary stabilization for a KdV system at noncritical lengths. Rapid stabilization was achieved through static and dynamic feedback controls constructed by means of Gramian operators, while time-dependent feedback laws were designed to drive the solutions to equilibrium in finite time. These results show that the stabilization theory for the KdV equation includes not only uniform exponential decay, but also arbitrarily prescribed exponential decay rates and finite-time stabilization under suitable feedback mechanisms.

A different Lyapunov-based approach was proposed by Baudouin, Cr\'epeau, and Valein~\cite{BaudouinCrepeauValein}, who investigated the stabilization of the nonlinear KdV equation with boundary time-delay feedback. They established exponential stability under suitable assumptions on the delay coefficients and the length of the spatial domain, obtaining an explicit estimate of the decay rate via a Lyapunov functional. They also developed an alternative proof based on an observability inequality, valid for any noncritical length, although without an explicit characterization of the decay rate.

Finally, in \cite{GaMu2025}, the authors investigated the boundary stabilization of the linearized KP-II equation posed on a bounded domain, without length restrictions on the domain. Their main objective was to design a suitable feedback law ensuring the existence and exponential stabilization of solutions in the energy space, without imposing restrictions on the length of the spatial domain. In addition, they analyzed the interplay between the drift term, the dispersive effects, and the boundary conditions, highlighting their role in the stabilization mechanism.

Despite this extensive literature, analogous controllability and stabilization results with one control are, to the best of our knowledge, still unavailable for the KP-II equation posed on a bounded rectangle with a single boundary control acting on one side of the domain. The main objective of the present work is to fill this gap by establishing exact boundary controllability and exponential stabilization, together with a characterization of the corresponding critical lengths.

\color{black}

\subsection{Notations and main results} Motivated by the previous discussion, we now present the main results of the paper.  

Let us introduce the functional space required for our analysis before presenting answers to our questions. Given $\Omega \subset \mathbb{R}^2$ a rectangular domain, let us denote the standard Hilbert space by $H^k(\Omega)$ endowed with its usual norm $\| \cdot \|_{H^k(\Omega)}$. 

Since the KP-II equation~\eqref{eq:KP} presents the non-local operator acting only in the longitudinal variable $x$, we introduce the anisotropic Hilbert space
$$
H_x^k(\Omega):=
	\left\lbrace%
	 	\varphi \colon \partial_x^j \varphi \in L^2(\Omega),\ \text{for } 0\leq j \leq k
	 \right\rbrace ,
$$
equipped with the norm $$\left\lVert \varphi \right\rVert_{H_x^k(\Omega)}^2 =\sum_{j=0}^k \left\lVert \partial_x^j \varphi\right\rVert_{L^2(\Omega)}^2.$$ Finally,  $H_{x0}^k(\Omega)$ will denote the closure of $C_0^\infty(\Omega)$ in $H_x^k(\Omega)$. 

For $T>0$, let us introduce the following set
\begin{equation*}
\mathcal{B}_{H^k} := C\left([0,T] ; L^2(\Omega)\right)\cap L^2(0,T ; H_{x0}^k(\Omega)),
\end{equation*}
endowed with its natural norm
\begin{equation*}
\| u \|_{\mathcal{B}_{H^k}} :=
		 \max_{t\in[0,T]} \| u(\cdot,\cdot, t) \|_{L^2(\Omega)} + \left(\int_0^T \| u(\cdot,\cdot,t)\|_{H_{x0}^k(\Omega)}^2\,dt\right)^{\frac{1}{2}}.
\end{equation*}

By following the ideas introduced by \cite{Rosier}, and writing
$$
		\mathfrak{a}:=m_1+2m_2+m_3,
$$
	consider the set
	\begin{equation}\label{critical-r}
		\mathcal{R}:=
		\left\lbrace
		\frac{\pi}{\sqrt{2}}
		\sqrt{\mathfrak{a}^2+2m_1^2+2m_3^2}
		\ ;\
		(m_1,m_2,m_3)\in\mathcal{M}
		\right\rbrace,
	\end{equation}
	where 
\begin{equation}\label{eq:cond-M}
\mathcal{M}
=
\left\{
(m_1,m_2,m_3)\in(\mathbb{N}^*)^3
\;\middle|\;
\begin{aligned}
&|m_1-m_3|>2m_2>0,\\[1mm]
&-\frac{
\left(\mathfrak{a}^2-4m_1^2\right)
\left(\mathfrak{a}^2-4m_3^2\right)
}{
8\left(\mathfrak{a}^2+2m_1^2+2m_3^2\right)
}
=n^2,
 \text{ for some } n\in\mathbb{N}^*
\end{aligned}
\right\}.
\end{equation}
		Note that the first condition in the set \eqref{eq:cond-M} guarantees that the second condition is positive, since a direct computation gives
	\begin{multline*}
		-\left(\mathfrak{a}^2-4m_1^2\right)\left(\mathfrak{a}^2-4m_3^2\right)
		\\
		=
		\left(3m_1+2m_2+m_3\right)
		\left((m_1-m_3)^2-4m_2^2\right)
		\left(m_1+2m_2+3m_3\right).
	\end{multline*}

The first main result answers Problem~A by establishing the exact boundary controllability of the KP-II equation.
\begin{theorem}\label{th:Control} Let $L\in(0,+\infty) \setminus \mathcal{R}$ and $T >0$. Then the KP-II equation~\eqref{eq:KP}-\eqref{eq:BC} is exactly controllable in time $T$, that  is, for any $u_0, u_T \in L^2(\Omega)$, there exists $h(y,t)\in L^2((0,T)\times(0,L))$ such that the mild solution $u$ of the KP-II equation~\eqref{eq:KP}-\eqref{eq:BC} satisfies $u(\cdot,\cdot,T) = u_T(x,y)$.
\end{theorem}

Next, by choosing a suitable feedback  damping mechanism as $h(y,t) = - \alpha u_x(0,y,t)$ with the constraint $0<|\alpha|\leq 1$, we get the boundary conditions~\eqref{eq:BC} with a feedback damping mechanism
\begin{equation}\label{eq:BC.f}
\begin{cases}
u(0,y,t) = u(L,y,t) = 0, & y\in(0,L) ,\ t\in(0,T),\\
u_x(L,y,t) = -\alpha u_x(0,y,t), & y\in(0,L) ,\ t\in(0,T),\\
u(x,0,t) = u(x,L,t) = 0, & x\in(0,L) ,\ t\in(0,T),\\
u(x,y,0) = u_0(x,y),& (x,y)\in\Omega.
\end{cases}
\end{equation}
Therefore, taking into account the boundary conditions~\eqref{eq:BC.f} we get that 
\begin{equation}\label{eq:S.En'}
\begin{aligned}
\frac{\mathrm{d}}{\mathrm{d} t} E(t) 
	& = - \frac{1-\alpha^2}{2} \int_0^L u_x^2(0,y)\,dy - \frac{1}{2}\int_0^L \left(\partial_x^{-1}u_y(0,y)\right)^2\,dy 
	 \leq 0.
\end{aligned}
\end{equation}
This allows us to establish exponential stabilization, which constitutes the second main result of this work and provides an answer to Problem B.
\begin{theorem}[Uniform exponential stabilization]\label{th:S.Exp}
Let $L\in(0,+\infty) \setminus \mathcal{R}$. Then, for any initial data $u_0\in L^2(\Omega)$  the energy $E(t)$, given by \eqref{eq:S.En}, associated with KP-II system \eqref{eq:KP}-\eqref{eq:BC.f}  decays exponentially. 
\end{theorem}

Theorem~\ref{th:S.Exp} guarantees the exponential decay of the energy through the
compactness--uniqueness method. However, this approach does not provide an explicit
decay rate and depends on the critical lengths. To avoid the restriction on some lengths partially and obtain quantitative information on the asymptotic behavior of the solutions, we next employ a Lyapunov approach. Therefore, our third main result is the following.

\begin{theorem}[Explicit decay rate]\label{th:S.Lyap}
Let $0<L<\sqrt{3}$, and let $\gamma>0$ satisfy
\begin{equation}\label{eq:S.Lyap.c2}
0<\gamma<\frac{1-\alpha^2}{L\alpha^2}.
\end{equation}
Set
\[
\rho:=1+L\gamma.
\]
Then, for every decay rate $\nu$ satisfying
\begin{equation}\label{eq:S.Lyap.c1}
0<\nu<
\frac{\gamma(3-L^2)}
{(1+L\gamma)L^2},
\end{equation}
and every initial data $u_0\in L^2(\Omega)$, the energy $E(t)$ defined in~\eqref{eq:S.En}, associated with the solution of the KP-II system~\eqref{eq:KP}–\eqref{eq:BC.f}, satisfies
\begin{equation}
E(t)\leq \rho E(0)e^{-\nu t},
\qquad t\geq0.
\end{equation}
In particular, the energy decays exponentially at any rate $\nu$ satisfying~\eqref{eq:S.Lyap.c1}.
\end{theorem}

The proofs of Theorems \ref{th:Control} and \ref{th:S.Exp} are a consequence of the compactness-uniqueness argument due to J.-L. Lions~\cite{Lions88}. The main idea is to prove an observability inequality that reduces our problem to proving a unique continuation property. This unique continuation is achieved thanks to the good properties of the spectral problem associated with both problems. 

Note that, unlike Theorem \ref{th:S.Exp}, Theorem \ref{th:S.Lyap} does not ensure the critical length phenomenon; however, a restrictive assumption on the length $L$ is necessary to control the energy of the system. This happens because the Lyapunov method does not involve the spectrum of the operator associated with the stabilization problem. We encourage the reader to see Section \ref{Sec3} for more details. 

The exact control property using the H.U.M. approach depends strongly on the \emph{a priori} estimates and the regularity results, and unfortunately, the focusing case, KP-I, does not allow the use of these tools to obtain equivalent results. Many authors have tried in the last decades to give a new set of critical lengths for the dispersive system; we can cite, for example,  \cite{CaChSoGo2023,Capistrano2019a,Cerpa-Zh,GG} and therein. After the pioneering work \cite{Rosier}, only the following articles were able to characterize the explicit critical length in the context of the bounded domain for KdV-type equations \cite{Capistrano16,Capistrano17,Capistrano2019,CaCaZh} and \cite{Doronin} for the 2D Zakharov-Kuznetsov equation posed on bounded rectangles and on a strip. So, few works were able to present a new and complete picture of the critical length set for dispersive-type systems in one- or two-dimensional cases.

Thus, our main novelty is to provide an explicit characterization of the set \eqref{critical-r} for the Kadomtsev-Petviashvili II equation. The critical lengths are described by a closed formula in the parameters $m_1,m_2,m_3$, subject to a single Diophantine compatibility condition which, for any prescribed bound on $L$, can be checked in finitely many steps.

\subsection{Paper's outline}
Our paper is organized into five sections, including the introduction. In Section \ref{Sec1}, we address the well-posedness problem associated with the KP-II system. Section \ref{Sec2} presents two types of observability inequalities, which allow us to explicitly characterize the critical length described by the set \eqref{critical-r}. These observability inequalities are crucial for establishing Theorems \ref{th:Control} and \ref{th:S.Exp}, which are proved in Section \ref{Sec3}. In the same section, we use Lyapunov's method to prove Theorem \ref{th:S.Lyap}, which relaxes the critical length condition for the stabilization problem. Finally, Section \ref{Sec4} discusses some properties of the set $\mathcal{R}$ and open issues.

\section{Well-posedness for the KP-II}\label{Sec1}

\subsection{Homogeneous system} Let us start with the well-posedness of the following problem 
\begin{equation}\label{eq:KP_h}
\begin{cases}
u_t + u_x + u_{xxx} +\partial_x^{-1}(u_{yy}) = 0, & (x,y)\in \Omega ,\ t\in(0,T), \\
u(0,y,t) = u(L,y,t) = u_x(L,y,t) = 0, & y\in(0,L) ,\ t\in(0,T)\\
u(x,0,t) = u(x,L,t) = 0, & x\in(0,L) ,\ t\in(0,T), \\
u(x,y,0) = u_0(x,y), & (x,y)\in \Omega.
\end{cases}
\end{equation}
Associated with this system, denote the operator $A\colon D(A) \subset L^2(\Omega) \to L^2(\Omega)$ defined by
$
Au = -u_x - u_{xxx} - \partial_x^{-1} u_{yy},$ with dense domain given by
\begin{equation}\label{eq:Adom}
D(A) := \left\lbrace u \in H_x^3(\Omega)  \cap H^2(\Omega)  \left\lvert
\begin{aligned}
\ u(0,y) = u(L,y) = u(x,0) = u(x,L) = u_x(L,y) = 0
\end{aligned}\right.\right\rbrace.
\end{equation}
We present some useful results to establish the linear system's well-posedness.
\begin{lemma}\label{le:Aadj}
The operator $A$ is closed and the adjoint operator $A^\ast\colon D(A^\ast)\subset L^2(\Omega)\to L^2(\Omega)$ is given by\footnote{Observe that for the adjoint problem, the operator $\left(\partial_x^{-1}\right)^{\ast}$ is defined as $
\left({\partial_x^{-1}}\right)^\ast \varphi(x,y) = \psi(x,y)$ such that $\psi(0,y) = 0\text{ and }\psi_x(x,y)=\varphi(x,y).$ For this case, the definition is equivalent to $\left(\partial_x^{-1}\right)^{\ast} \varphi(x,y,t) = \int_0^x \varphi(s,y,t)\,ds$.}
\begin{equation}\label{eq:Aadj}
A^\ast \theta =
		\theta_x + \theta_{xxx}+ \left(\partial_x^{-1}\right)^{\ast} \theta_{yy}
\end{equation}
with dense domain
\begin{align}\label{eq:Aadom}
D(A^\ast):=
&\left\lbrace \theta \in H_x^3(\Omega)  \cap H^2(\Omega)  \left\lvert
\begin{aligned}
\theta(0,y) = \theta(L,y) = \theta(x,0) = \theta(x,L) = \theta_x(0,y) = 0
\end{aligned}\right.  \right\rbrace. 
\end{align}
\end{lemma}
\begin{proof}
Let $u\in D(A)$ and $\theta \in D(A^\ast)$. We consider two functions, $f, g$ such that $u_y(x,y)=f_x(x,y)$ and $\theta_y(x,y)=g_x(x,y)$, with $f(L,y)=0$ and $g(0,y)=0$.  Firstly, we estimate the product of $\theta$ by the nonlocal term ${\partial_x^{-1}}u_{yy}$, so we have, by integration by parts and the trace properties, that
\begin{equation*}
\begin{aligned}
-\int_{\Omega} \theta(x,y){\partial_x^{-1}}u_{yy}(x,y)\,dx\,dy 
 &= -\int_{\Omega}  \theta(x,y)f_y(x,y)\,dx\,dy \\
 &=	\int_{\Omega}  u(x,y)\left({\partial_x^{-1}}\right)^\ast  \theta_{yy}(x,y)\,dx\,dy.%
\end{aligned}
\end{equation*}
Consequently, we can estimate the duality product $\left\langle
 		A u, \theta
 	\right\rangle_{L^2(\Omega)}$ as follows 
\begin{equation*}
\begin{aligned}
	\left\langle
 		A u, \theta
 	\right\rangle_{L^2(\Omega)}  
		& = \int_{\Omega} u(x,y)\left(\theta_x(x,y) + \theta_{xxx}(x,y) + {\partial_x^{-1}}\theta_{yy}(x,y)\right)\,dx\,dy \\
		& =\left\langle%
 		u, A^\ast \theta
 	\right\rangle_{L^2(\Omega)}.
\end{aligned}
\end{equation*}
Moreover, note that $A^{\ast\ast}=A$, then $A$ is closed.

Finally, let us highlight that $D(A)$ is dense in $L^2(\Omega)$. Indeed, since $C_0^\infty(\Omega)\subset D(A)$ and $C_0^\infty(\Omega)$ is dense in $L^2(\Omega)$, it follows immediately that $\overline{D(A)}^{L^2(\Omega)}=L^2(\Omega)$. More precisely, to verify the inclusion, let $\varphi\in C_0^\infty(\Omega)$. Then $\varphi\in H_x^3(\Omega)\cap H^2(\Omega)$ and all the boundary conditions in the definition of $D(A)$ are satisfied, namely
$$
\varphi(0,y)=\varphi(L,y)=\varphi(x,0)=\varphi(x,L)=0,
\qquad
\varphi_x(L,y)=0.
$$
Moreover, the nonlocal term $\partial_x^{-1}\varphi_{yy}$ belongs to $L^2(\Omega)$.
Therefore, $\varphi\in D(A)$, and hence $D(A)$ is dense in $L^2(\Omega)$.\color{black}
\end{proof}

Let us now show that $A$ generates a $C_0$-semigroup in $L^2(\Omega)$.

\begin{proposition}\label{pr:Sem_Dis}
The operator $A$ is the infinitesimal generator of a $C_0$-semigroup in $L^2(\Omega)$.
\end{proposition}

\begin{proof}
Let $u\in D(A)$ and call $f_x(x,y)=u_y(x,y)$ with $f(L,y)=0$. From~\eqref{eq:Adom} it follows, thanks to the integration by parts, that
\begin{equation*}
\begin{aligned}
\left\langle Au,u	 \right\rangle_{L^2(\Omega)}
	&=\frac{1}{2}\int_0^L [u_x^2 (L,y)- u_x^2(0,y)]\,dy 		+\frac{1}{2}\int_0^L[f^2(L,y) - f^2(0,y)]\,dy. 
\end{aligned}
\end{equation*}
Therefore, we have that
\begin{equation}
\label{eq:Adis}
\begin{aligned}
	\left\langle%
		Au,u%
	\right\rangle_{L^2(\Omega)}
				&= - \frac{1}{2}\int_0^L u_x^2(0,y)dy 
				-\frac{1}{2}\int_0^L(\partial_x^{-1} u_y(0,y))^2 	\,dy   \leq 0.
\end{aligned}
\end{equation}
Similarly, let $\theta\in D(A^\ast)$ and $g_x(x,y)=\theta_y(x,y)$ with $g(0,y)=0$. Then
\begin{equation*}
\begin{aligned}
\left\langle \theta, A^\ast \theta \right\rangle_{L^2(\Omega)}
	&= -\frac{1}{2}\int_0^L [\theta_x^2 (L,y)- \theta_x^2(0,y)]\,dy - \frac{1 }{2}\int_{\Omega} (g^2(x,y))_x \,dx\,dy .
\end{aligned}
\end{equation*}
From~\eqref{eq:Aadom}, it yields that
\begin{equation*}
\left\langle \theta, A^\ast \theta \right\rangle_{L^2(\Omega)} 
	= - \frac{1}{2}\int_0^L \theta_x^2(L,y)\,dy -\frac{1}{2}\int_0^L\left((\partial_x^{-1})^{\ast} \theta_y(L,y)\right)^2 \,dy\leq 0.  
\end{equation*}
Thus, from  \cite[Corollary 4.4, page 15]{Pazy}, the proposition holds. 
\end{proof}

Now, we can state the following theorem regarding the existence of solutions to the Cauchy problem.
\begin{theorem}
For each initial data $u_0\in L^2(\Omega)$, there exists a unique mild solution $u\in C\left([0,\infty), L^2(\Omega)\right)$ for the system~\eqref{eq:KP_h}. Moreover, if the initial data $u_0\in D(A)$ the solution $u$ belongs to  
$ C\left([0,\infty); D(A)\right)\cap C^1\left([0,\infty); L^2(\Omega)\right).$
\end{theorem}
\begin{proof}
From Proposition~\ref{pr:Sem_Dis}, it follows that $A$ generates a strongly continuous semigroup of
contractions $\{S(t)\}_{t\geq 0}$ in $L^2(\Omega)$ (see  \cite[Corollary 1.4.4]{Pazy}). 
\end{proof}

The next proposition provides valuable estimates for the mild solutions of the equation \eqref{eq:KP_h}, including the energy estimate, Kato smoothing effect, and the existence of the traces.
\begin{proposition}\label{pr:Kato}
Let $u_0 \in L^2(\Omega)$, then the map $u_0 \in L^2(\Omega)\mapsto u \in\mathcal{B}_{H^1}$ is well-defined, continuous and the corresponding solution $u$ of system~\eqref{eq:KP_h} satisfies
\begin{equation}\label{eq:reg-aa}
\| u(\cdot,\cdot,t)\|_{L^2(\Omega)} \leq \| u_0 \|_{L^2(\Omega)}, \quad\forall t\in[0,T]
\end{equation}
and
\begin{equation}\label{eq:reg}
\| u \|_{L^2(0,T; H_x^1(\Omega))}^2 + \| {\partial_x^{-1}}u_y\|_{L^2(0,T; L^2(\Omega))}^2 \leq C(T,L) \| u_0\|_{L^2(\Omega)}^2,
\end{equation}
for some positive constant $C:= C(T,L) = (2T+L)>0$. 
Moreover, if the initial data $u_0\in L^2(\Omega)$, we get the following trace estimate
\begin{equation}\label{eq:traces}
\|u_x(0,\cdot,\cdot)\|_{L^2((0,T)\times (0,L))}^2  \leq  \| u_0 \|_{L^2(\Omega)}^2
\end{equation}
and the estimate
\begin{equation}\label{eq:Time_est}
\begin{aligned}
\int_{\Omega} u_0^2(x,y)\,dx\,dy 
	& \leq \frac{1}{T} \int_0^T\int_{\Omega} u^2\,dx\,dy\,dt + \int_0^T\int_0^L u_x^2(0,y,t)\,dy\,dt\\
	& \quad + \int_0^T\int_0^L \left({\partial_x^{-1}}u_y(0,y,t)\right)^2\,dy\,dt.
\end{aligned}
\end{equation}
\end{proposition}

\begin{proof}

Observe that~\eqref{eq:Adis} implies
\begin{equation*}
\begin{aligned}
\frac{\mathrm{d}}{\mathrm{d} t} \frac{1}{2}\int_{\Omega} u^2\,dx\,dy 
	  & =\int_{\Omega} u(t) \left(-u_{x}(t) - u_{xxx}(t)- {\partial_x^{-1}} u_{yy}(t)  \right)  \,dx\,dy  \\
	 & = \left\langle Au(t),u(t) \right\rangle_{L^2(\Omega)} \leq 0.
\end{aligned}
\end{equation*}
Integrating over $[0,s]$, for $0\leq s \leq T$, immediately yields \eqref{eq:reg-aa}.

Next, we prove the Kato smoothing estimate, namely, $u\in L^2(0,T;H_x^1(\Omega)).$ Multiplying the equation in~\eqref{eq:KP_h} by $xu$ and integrating by parts over $\Omega\times(0,T)$, we obtain
\begin{equation*}
\begin{aligned}
\frac{3}{2}\int_0^T\int_{\Omega}u_x^2\,dx\,dy\,dt
+\frac{1}{2}\int_0^T\int_{\Omega}
\left(\partial_x^{-1}u_y\right)^2\,dx\,dy\,dt
&=
\frac{1}{2}\int_0^T\int_{\Omega}u^2\,dx\,dy\,dt\\
&\quad
+\frac{1}{2}\int_{\Omega}xu_0^2(x,y)\,dx\,dy\\
&\quad
-\frac{1}{2}\int_{\Omega}xu^2(x,y,T)\,dx\,dy.
\end{aligned}
\end{equation*}

Since
\begin{equation*}
\begin{split}
\frac{3}{2}\int_0^T\int_{\Omega}u_x^2\,dx\,dy\,dt
&+\frac{1}{2}\int_0^T\int_{\Omega}
\left(\partial_x^{-1}u_y\right)^2\,dx\,dy\,dt\\
&\geq
\frac{1}{2}
\left(
\|u_x\|_{L^2(0,T;L^2(\Omega))}^2
+
\|\partial_x^{-1}u_y\|_{L^2(0,T;L^2(\Omega))}^2
\right),
\end{split}
\end{equation*}
we deduce that
\begin{equation*}
\begin{aligned}
\frac{1}{2}
\left(
\|u_x\|_{L^2(0,T;L^2(\Omega))}^2
+
\|\partial_x^{-1}u_y\|_{L^2(0,T;L^2(\Omega))}^2
\right)
&\leq
\frac{1}{2}\int_0^T\int_{\Omega}u^2\,dx\,dy\,dt\\
&\quad
+\frac{1}{2}\int_{\Omega}xu_0^2(x,y)\,dx\,dy\\
&\quad
-\frac{1}{2}\int_{\Omega}xu^2(x,y,T)\,dx\,dy.
\end{aligned}
\end{equation*}
The last term on the right-hand side is nonpositive and may therefore be discarded.

Moreover, by~\eqref{eq:reg-aa} and the fact that $0\leq x\leq L$, we have
\begin{equation}\label{eq:kato_aux}
\int_0^T\int_\Omega u^2\,dx\,dy\,dt
\leq
T\|u_0\|_{L^2(\Omega)}^2,
\qquad
\int_\Omega xu_0^2\,dx\,dy
\leq
L\|u_0\|_{L^2(\Omega)}^2.
\end{equation}
Hence,
\begin{equation*}
\frac12
\left(
\|u_x\|_{L^2(0,T;L^2(\Omega))}^2
+
\|\partial_x^{-1}u_y\|_{L^2(0,T;L^2(\Omega))}^2
\right)
\leq
\frac{T+L}{2}\|u_0\|_{L^2(\Omega)}^2,
\end{equation*}
and therefore
\begin{equation*}
\|u_x\|_{L^2(0,T;L^2(\Omega))}^2
+
\|\partial_x^{-1}u_y\|_{L^2(0,T;L^2(\Omega))}^2
\leq
(T+L)\|u_0\|_{L^2(\Omega)}^2.
\end{equation*}

Consequently,
\begin{equation*}
\begin{aligned}
\|u\|_{L^2(0,T;H_x^1(\Omega))}^2
&=
\|u\|_{L^2(0,T;L^2(\Omega))}^2
+
\|u_x\|_{L^2(0,T;L^2(\Omega))}^2\\
&\leq
T\|u_0\|_{L^2(\Omega)}^2
+
\|u_x\|_{L^2(0,T;L^2(\Omega))}^2.
\end{aligned}
\end{equation*}
Combining the last two estimates yields
\begin{equation*}
\|u\|_{L^2(0,T;H_x^1(\Omega))}^2
+
\|\partial_x^{-1}u_y\|_{L^2(0,T;L^2(\Omega))}^2
\leq
(2T+L)\|u_0\|_{L^2(\Omega)}^2,
\end{equation*}
which proves~\eqref{eq:reg}.

For the trace estimate, we multiply the equation \eqref{eq:KP_h} by $u$ and integrate by parts over $ \Omega \times (0,T)$,
\begin{multline*}
\int_0^T \int_0^L u_x^2(0,y,t)\,dy\,dt + \int_0^T\int_0^L \left({\partial_x^{-1}}u_y(0,y,t)\right)^2\,dy\,dt \\= \int_{\Omega} u_0^2(x,y)\,dx\,dy - \int_{\Omega} u^2(x,y,T)\,dx\,dy.
\end{multline*}
Since the second term on the left-hand side is nonnegative, we conclude that
\begin{equation*}
\begin{aligned}
\|u_x(0,\cdot,\cdot)\|_{L^2((0,T)\times (0,L))}^2 
	& \leq
	\|u_x(0,\cdot,\cdot)\|_{L^2((0,T)\times (0,L))}^2 + \left\|{\partial_x^{-1}}u_y(0,\cdot,\cdot)\right\|_{L^2((0,T)\times (0,L))}^2 \\
	 & \leq
		 \| u_0 \|_{L^2(\Omega)}^2,
\end{aligned}
\end{equation*}
and \eqref{eq:traces} holds.

Lastly, multiplying the equation \eqref{eq:KP_h} by $(T-t)u$ and integrating by parts in $ \Omega \times (0,T)$ yields that
$$
\begin{aligned}
\int_{\Omega} u_0^2(x,y)\,dx\,dy 
	& \leq \frac{1}{T} \int_0^T\int_{\Omega} u^2\,dx\,dy\,dt + \int_0^T\int_0^L u_x^2(0,y,t)\,dy\,dt\\
	& \quad + \int_0^T\int_0^L \left({\partial_x^{-1}}u_y(0,y,t)\right)^2\,dy\,dt,
\end{aligned}
$$
showing \eqref{eq:Time_est}, and the result is proven.
\end{proof}

\subsection{Adjoint system}
We consider the following time-backward homogeneous problem
\begin{equation}\label{eq:tbk_h.1}
\begin{cases}
\theta_t + \theta_x +\theta_{xxx} + \left({\partial_x^{-1}}\right)^\ast \theta_{yy} = 0 , & (x,y)\in \Omega,\ t\in (0,T),\\
\theta(0,y,t) = \theta(L,y,t) = \theta_x(0,y,t) = 0, & y\in (0,L) ,\ t\in (0,T), \\
\theta(x,0,t) = \theta(x,L,t) = 0, & x\in (0,L) ,\ t\in (0,T), \\
\theta(x,y,T) = \theta_T(x,y), & (x,y)\in \Omega .
\end{cases}
\end{equation}
Note that, by the variable change $t\mapsto T-t$ the problem~\eqref{eq:tbk_h.1} can be formulated as
\begin{equation}\label{eq:tbk_h.2}
\begin{cases}
\theta_t - \theta_x - \theta_{xxx} - \left({\partial_x^{-1}}\right)^\ast \theta_{yy} = 0 , & (x,y)\in \Omega,\ t\in (0,T),\\
\theta(0,y,t) = \theta(L,y,t) = \theta_x(0,y,t) = 0, & y\in (0,L) ,\ t\in (0,T), \\
\theta(x,0,t) = \theta(x,L,t) = 0, & x\in (0,L) ,\ t\in (0,T), \\
\theta(x,y,0) = \theta_0(x,y), & (x,y)\in \Omega .
\end{cases}
\end{equation}
That is equivalent to abstract Cauchy problem $\theta_t = A^\ast\theta$ with initial data $\theta(x,y,0) = \theta_0(x,y)$, where $A^\ast$ is defined by~\eqref{eq:Aadj} with dense domain~\eqref{eq:Aadom}, by using Lemma~\ref{le:Aadj} we can establish the well-posedness for the problem~\eqref{eq:tbk_h.2} and consequently the well-posedness for~\eqref{eq:tbk_h.1} yields. Precisely, we get the following result.
\begin{theorem}
For each initial data $\theta_0\in L^2(\Omega)$, there exists a unique mild solution $\theta\in C\left([0,\infty), L^2(\Omega)\right)$ for the problem~\eqref{eq:tbk_h.2}. Moreover, if the initial data $\theta_0\in D(A^\ast)$ the solutions are classical such that
$\theta\in C\left([0,\infty), D(A^\ast)\right)\cap C^1\left([0,\infty), L^2(\Omega)\right).$
\end{theorem}
\subsection{KP-II equation with feedback}
Let us first pay attention to the existence and regularity issues for the solutions to~\eqref{eq:KP} with a feedback mechanism given by~\eqref{eq:BC.f}. Gathering this information, we get the system
\begin{equation}\label{eq:S.KP}
\begin{cases}
u_t + u_x + u_{xxx} +\partial_x^{-1}(u_{yy}) = 0, & (x,y)\in \Omega ,\ t\in(0,T), \\
u(0,y,t) = u(L,y,t) = 0, & y\in(0,L) ,\ t\in(0,T)\\
u(x,0,t) = u(x,L,t) = 0,\quad u_x(L,y,t) =  -\alpha u_x(0,y,t), &  x\in(0,L) ,\ t\in(0,T)\\
u(x,y,0) = u_0(x,y), & (x,y)\in \Omega.
\end{cases}
\end{equation}
We introduce the operator $\tilde{A}: D(\tilde{ A})\subset L^2(\Omega) \to L^2(\Omega)$ defined by
$
\tilde{A} u : = -u_x - u_{xxx} - {\partial_x^{-1}}u_{yy}
$
on the dense domain given by
$$
D(\tilde{A}) := \left\lbrace u \in H_x^3(\Omega)  \cap H^2(\Omega)  \left\lvert
\begin{aligned}
&\ u(0,y) = u(L,y) = u(x,0) = u(x,L) = 0\\
&\ u_x(L,y) = -\alpha u_x(0,y), \ 0 < |\alpha| < 1
\end{aligned}\right.  \right\rbrace.
$$
By analogous arguments as in Lemma~\ref{le:Aadj} we get that $\tilde{A}$ is closed and the adjoint  $\tilde{A}^\ast: D(\tilde{A}^\ast) \subset L^2(\Omega) \to L^2(\Omega)$ by
$
\tilde{A}^\ast v = v_x + v_{xxx} + \left(\partial_x^{-1}\right)^\ast v_{yy}
$
with dense domain
$$
D(\tilde{A}^\ast) := \left\lbrace v \in H_x^3(\Omega)  \cap H^2(\Omega)  \left\lvert
\begin{aligned}
&\ v(0,y) = v(L,y) = v(x,0) = v(x,L) = 0\\
&\ v_x(0,y) = - \alpha v_x(L,y), \ 0 < |\alpha| < 1
\end{aligned}\right.  \right\rbrace.
$$

{
Notice that, following the same reasoning used to prove the density of $D(A)$ in Lemma~\ref{le:Aadj}, we also obtain that $D(\widetilde A)$ is dense in $L^2(\Omega)$.}

\begin{proposition}\label{pr:S.diss}
The operator $\tilde{A}$ is the infinitesimal generator of a $C_0$-semigroup in $L^2(\Omega)$.
\end{proposition}
\begin{proof}
Similarly to the proof of the Proposition \ref{pr:Sem_Dis}, we consider  $u\in D(\tilde{A})$ and then
\begin{align*}
\left\langle \tilde{A}u,u	\right\rangle_{L^2(\Omega)}
			=- \frac{(1-\alpha^2)}{2}\int_0^L u_x^2(0,y)\,dy  -\frac{1}{2}\int_0^L(\partial_x^{-1} u_y(0,y))^2 \,dy \leq 0 . 
\end{align*}
Similarly, let $v\in D(\tilde{A}^\ast)$ therefore
\begin{align*}
\left\langle v, \tilde{A}^\ast v \right\rangle_{L^2(\Omega)}
			= - \frac{(1-\alpha^2)}{2}\int_0^L v_x^2 (L,y)\,dy  - \frac{1}{2} \int_0^L \left(({\partial_x^{-1}})^\ast v_y(L,y)\right)^2\,dy \leq 0.
\end{align*}
Thus, from  \cite[Corollary 4.4, page 15]{Pazy}, the proposition holds. 
\end{proof}

Therefore, we are in a position to present the subsequent theorem concerning the existence of solutions for the Cauchy abstract problem:
\begin{theorem}
For the initial data $u_0\in L^2(\Omega)$, there exists a unique solution $u\in C\left([0,\infty), L^2(\Omega)\right)$ for the KP-II equation~\eqref{eq:S.KP}. Moreover, if the initial data $u_0\in D(\tilde{A})$ the solutions are classical such that $u\in C\left([0,\infty), D(\tilde{A})\right)\cap C^1\left([0,\infty), L^2(\Omega)\right).$
\end{theorem}

\begin{proof}
From Proposition~\ref{pr:S.diss}, it follows that $\tilde{A}$ generates a strongly continuous semigroup of
contractions $(\tilde{S}(t))_{t\geq 0}$ in $L^2(\Omega)$ (see \cite[Corollary 1.4.4]{Pazy}). 
\end{proof}

\begin{proposition}
For any mild solution of~\eqref{eq:S.KP}, the energy $E(t)$ is non-increasing and there exists a constant $C>0$ such that
\begin{equation}\label{eq:S.E'(t)}
\begin{aligned}
\frac{\mathrm{d}}{\mathrm{d} t}E(t) &\leq -C
	\left[%
		\int_0^L u_{x}^2(0,y,t)\,dy
		+\int_0^L({\partial_x^{-1}} u_y(0,y,t))^2\,dy 
	\right]
\end{aligned}
\end{equation}
where $ C=C(\alpha)$ is given by $C=\min\left\lbrace \frac{1-\alpha^2}{2},\frac{1}{2}\right\rbrace$.
\end{proposition}
\begin{proof}
Note that 
\begin{align*}
\frac{\mathrm{d}}{\mathrm{d} t} E(t) & =\int_{\Omega} u(t) \left( - u_x(t) - u_{xxx}(t)- {\partial_x^{-1}} u_{yy}(t)  \right)  dxdy= \left\langle \tilde{A}u(t),u(t)	\right\rangle_{L^2(\Omega)}. 
\end{align*}
 From the proof of Proposition \ref{pr:S.diss}, we conclude~\eqref{eq:S.E'(t)}.
\end{proof}

The following estimates for the mild solution $u(\cdot) = \tilde{S}(\cdot)(u_0)$ of the KP-II equation provide useful information for the regularity and the asymptotic behavior. The first is standard energy estimates, and the last one reveals that a slightly different Kato smoothing effect holds.
\begin{proposition}
Let $u_0 \in L^2(\Omega)$ and the mild solution given by $u(\cdot) = \tilde{S}(\cdot)(u_0)$. Then,
\begin{equation}\label{eq:S.traces}
\begin{split}
\frac{1}{2}\int_{\Omega}  u_0^2(x,y) \,dx\,dy - \frac{1}{2}\int_{\Omega}  u^2(x,y,T) \,dx\,dy 
	=& \frac{1-\alpha^2}{2} \int_0^T\int_0^L u_x^2(0,y,t)\,dy\,dt \\&+ \frac{1}{2} \int_0^T\int_0^L  \left({\partial_x^{-1}} u_y(0,y,t)\right)^2 \,dy\,dt,
	\end{split}
\end{equation}
and
\begin{equation}\label{eq:S.traces_1}
\begin{split}
\frac{T}{2}\int_{\Omega} u_0^2(x,y)\,dx\,dy  =& \frac{1}{2} \int_0^T \int_{\Omega} u^2\,dx\,dy\,dt +\frac{1-\alpha^2}{2}\int_0^T\int_0^L (T-t) u_x^2(0,y,t)\,dy\,dt \\&+ \frac{1}{2} \int_0^T\int_0^L (T-t) \left({\partial_x^{-1}} u_y(0,y,t)\right)^2 \,dy\,dt,
	\end{split}
\end{equation}
 for any $T>0$. Moreover, we have
\begin{equation*}
\| u \|_{L^2(0,T; H_x^1(\Omega))}^2 \leq C(T,L,\alpha) \| u_0\|_{L^2(\Omega)}^2,
\end{equation*}
where $C:=C(T,L,\alpha) = \left(
2T+\frac{L}{1-\alpha^2}
\right)> 0$.
\end{proposition}
\begin{proof}
The identities \eqref{eq:S.traces} and \eqref{eq:S.traces_1} are obtained exactly as in
Proposition~\ref{pr:Kato}. More precisely, multiplying \eqref{eq:S.KP} by $u$ and integrating by parts yields \eqref{eq:S.traces}, whereas multiplying by $(T-t)u$ gives \eqref{eq:S.traces_1}\footnote{Notice that the only difference with Proposition~\ref{pr:Kato} is the contribution of the feedback law given boundary condition $ u_x(L,y,t)=-\alpha u_x(0,y,t),$ which produces the factor $(1-\alpha^2)$ in the boundary term.}.

Now, multiplying~\eqref{eq:S.KP} by $xu$ and integrating by parts over
$\Omega\times(0,T)$, we obtain
\begin{equation}\label{eq:S.Kato.1}
\begin{aligned}
&\frac{3}{2}\int_0^T\int_{\Omega}u_x^2\,dx\,dy\,dt
+\frac{1}{2}\int_0^T\int_{\Omega}
\left(\partial_x^{-1}u_y\right)^2\,dx\,dy\,dt\\
&\qquad =
\frac{1}{2}\int_{\Omega}xu_0^2(x,y)\,dx\,dy
-\frac{1}{2}\int_{\Omega}xu^2(x,y,T)\,dx\,dy\\
&\qquad\quad
+\frac{L\alpha^2}{2}
\int_0^T\int_0^L u_x^2(0,y,t)\,dy\,dt
+\frac{1}{2}\int_0^T\int_{\Omega}u^2\,dx\,dy\,dt.
\end{aligned}
\end{equation}

On the other hand, it follows from~\eqref{eq:S.traces} that
\begin{equation*}
\begin{aligned}
(1-\alpha^2)
\int_0^T\int_0^L u_x^2(0,y,t)\,dy\,dt
&\leq
(1-\alpha^2)
\int_0^T\int_0^L u_x^2(0,y,t)\,dy\,dt\\
&\quad
+\int_0^T\int_0^L
\left(\partial_x^{-1}u_y(0,y,t)\right)^2\,dy\,dt\\
&=
2E(0)-2E(T)\\
&\leq
2E(0)
=
\|u_0\|_{L^2(\Omega)}^2.
\end{aligned}
\end{equation*}
Hence,
\begin{equation}\label{eq:S.boundary.estimate}
\int_0^T\int_0^L u_x^2(0,y,t)\,dy\,dt
\leq
\frac{1}{1-\alpha^2}
\|u_0\|_{L^2(\Omega)}^2.
\end{equation}

Moreover, estimating from below the coefficients on the left-hand side
of~\eqref{eq:S.Kato.1}, we deduce that
\begin{equation*}
\begin{aligned}
\frac{1}{2}
\bigg(
\|u_x\|_{L^2(0,T;L^2(\Omega))}^2
+
\|\partial_x^{-1}u_y\|_{L^2(0,T;L^2(\Omega))}^2
\bigg)
&\leq
\frac{1}{2}\int_\Omega xu_0^2\,dx\,dy\\
&\quad
+\frac{L\alpha^2}{2}
\int_0^T\int_0^L u_x^2(0,y,t)\,dy\,dt\\
&\quad
+\frac{1}{2}\int_0^T\int_\Omega u^2\,dx\,dy\,dt,
\end{aligned}
\end{equation*}
where the nonpositive term
\[
-\frac{1}{2}\int_{\Omega}xu^2(x,y,T)\,dx\,dy
\]
has been discarded. Employing~\eqref{eq:kato_aux} and
\eqref{eq:S.boundary.estimate}, we obtain
\begin{equation*}
\begin{aligned}
\bigg(
\|u_x\|_{L^2(0,T;L^2(\Omega))}^2
+
\|\partial_x^{-1}u_y\|_{L^2(0,T;L^2(\Omega))}^2
\bigg)
&\leq
\left(
T
+L
+\frac{L\alpha^2}{(1-\alpha^2)}
\right)
\|u_0\|_{L^2(\Omega)}^2\\
&=
\frac{T(1-\alpha^2)+L}
{(1-\alpha^2)}
\|u_0\|_{L^2(\Omega)}^2.
\end{aligned}
\end{equation*}
Therefore,
\begin{equation}\label{eq:S.Kato.derivative}
\|u_x\|_{L^2(0,T;L^2(\Omega))}^2
+
\|\partial_x^{-1}u_y\|_{L^2(0,T;L^2(\Omega))}^2
\leq
\left(
T+\frac{L}{1-\alpha^2}
\right)
\|u_0\|_{L^2(\Omega)}^2.
\end{equation}

Finally, since
\[
\|u\|_{H_x^1(\Omega)}^2
=
\|u\|_{L^2(\Omega)}^2
+
\|u_x\|_{L^2(\Omega)}^2,
\]
we have
\begin{equation*}
\begin{aligned}
\|u\|_{L^2(0,T;H_x^1(\Omega))}^2
+
\|\partial_x^{-1}u_y\|_{L^2(0,T;L^2(\Omega))}^2
&=
\|u\|_{L^2(0,T;L^2(\Omega))}^2\\
&\quad
+\|u_x\|_{L^2(0,T;L^2(\Omega))}^2\\
&\quad
+\|\partial_x^{-1}u_y\|_{L^2(0,T;L^2(\Omega))}^2\\
&\leq
T\|u_0\|_{L^2(\Omega)}^2
+
\left(
T+\frac{L}{1-\alpha^2}
\right)
\|u_0\|_{L^2(\Omega)}^2\\
&=
\left(
2T+\frac{L}{1-\alpha^2}
\right)
\|u_0\|_{L^2(\Omega)}^2.
\end{aligned}
\end{equation*}
This proves the desired estimate.
\end{proof}

\subsection{Nonhomogeneous system} 
Let us consider problem~\eqref{eq:KP} with the nonhomogeneous boundary condition~\eqref{eq:BC} and initial data~$u_0(x,y)$. Proposition~\ref{pr:Kato} provides Kato smoothing and boundary trace estimates for the corresponding homogeneous problem. However, unlike the one-dimensional KdV equation, an analogue of the hidden regularity theory that underlies the boundary forcing operator approach is not currently available for the KP-II equation in the present setting\footnote{For an exposition of the boundary forcing operator approach and its connection with hidden regularity for the KdV equation, we refer the reader to~\cite{Bona03}.}. Such a theory typically provides additional time regularity for the boundary traces and plays a fundamental role in the treatment of nonhomogeneous boundary value problems with rough boundary data. Since this approach is not presently available, we instead establish the well-posedness of problem~\eqref{eq:KP}--\eqref{eq:BC} by means of the transposition (duality) method, which only requires the trace estimate obtained in Proposition~\ref{pr:Kato}.

As usual, we begin with a formal computation. Fix $\tau\in[0,T]$ and
$\theta_\tau\in D(A^*)$, where $D(A^*)$ is given by~\eqref{eq:Aadom},
and let $\theta$ be the solution of the backward adjoint problem
\eqref{eq:tbk_h.1} on the time interval $(0,\tau)$, with final time data
\[
\theta(x,y,\tau)=\theta_\tau(x,y).
\]
Multiplying~\eqref{eq:KP} by $\theta$, integrating over
$\Omega\times(0,\tau)$, and using the boundary conditions
\eqref{eq:BC} and those of the adjoint system, we formally obtain
\begin{equation}\label{eq:formal_transposition}
\begin{aligned}
0
&=
\int_\Omega u(x,y,\tau)\theta_\tau(x,y)\,dx\,dy
-
\int_\Omega u_0(x,y)\theta(x,y,0)\,dx\,dy\\
&\quad
-
\int_0^\tau\int_0^L
h(y,t)\theta_x(L,y,t)\,dy\,dt.
\end{aligned}
\end{equation}
Accordingly, for each $\tau\in[0,T]$, we define the linear functional
$L_\tau\colon D(A^*)\to\mathbb{R}$ by
\begin{equation}\label{eq:fc.tr}
L_\tau(\theta_\tau)
:=
\left\langle
u_0,\theta(0)
\right\rangle_{D(A^*)',D(A^*)}
+
\int_0^\tau\int_0^L
h(y,t)\theta_x(L,y,t)\,dy\,dt,
\end{equation}
where $\theta$ denotes the solution of the adjoint problem associated
with the final time data $\theta_\tau$.

Thus, identity~\eqref{eq:formal_transposition} can be rewritten as
\begin{equation}\label{transposition}
\left\langle
u(\tau),\theta_\tau
\right\rangle_{D(A^*)',D(A^*)}
=
L_\tau(\theta_\tau),
\qquad
\forall\,\theta_\tau\in D(A^*).
\end{equation}
This leads to the following definition.
\begin{definition}\label{transp-R}
Let $T>0$, $u_0\in D(A^*)'$, and
$h\in L^2((0,T)\times(0,L))$. We say that
\[
u\in C([0,T];D(A^*)')
\]
is a solution by transposition of~\eqref{eq:KP}--\eqref{eq:BC} if, for
every $\tau\in[0,T]$ and every $\theta_\tau\in D(A^*)$, identity
\eqref{transposition} holds, where $\theta$ is the solution of the
backward adjoint problem on $(0,\tau)$ with terminal datum
$\theta(\tau)=\theta_\tau$.
\end{definition}

We are now in a position to establish the well-posedness of the nonhomogeneous problem. Although the argument follows the same general strategy as that of~\cite[Proposition 2.6]{Capistrano2019}, we include the proof for completeness, since several modifications are required in the present setting. In particular, the proof relies only on the semigroup formulation of the adjoint problem together with the trace estimate established in Proposition~\ref{pr:Kato}, and does not require any hidden regularity properties.
\begin{proposition}
 Let  $T>0$, $u_0 \in D(A^*)'$ and $h(y,t)\in L^2((0,T)\times(0,L))$. Then there exists a unique solution $u\in C([0,T], D(A^*)')$ of \eqref{eq:KP}-\eqref{eq:BC} such that
\begin{equation}\label{eq:nh_est}
\| u \|_{L^{\infty}(0,T; D(A^*)')} \leq C \left(\|u_0 \|_{D(A^*)'} + \| h \|_{L^2((0,T)\times(0,L))}\right).
\end{equation}
\end{proposition}
\begin{proof}
Fix $\tau\in[0,T]$. We prove that the linear functional
$L_\tau$, defined in~\eqref{eq:fc.tr}, is continuous on $D(A^*)$,
endowed with the graph norm
\[
\|\varphi\|_{D(A^*)}
:=
\|\varphi\|_{L^2(\Omega)}
+
\|A^*\varphi\|_{L^2(\Omega)}.
\]
Let $\theta_\tau\in D(A^*)$, and denote by $\theta$ the corresponding
solution of the adjoint problem. If $(S^*(t))_{t\geq0}$ denotes the
$C_0$-semigroup generated by $A^*$ on $L^2(\Omega)$, then
\[
\theta(t)=S^*(\tau-t)\theta_\tau,
\qquad 0\leq t\leq\tau.
\]
Since $D(A^*)$ is invariant under $S^*(t)$ and
\[
A^*S^*(t)\varphi=S^*(t)A^*\varphi,
\qquad \varphi\in D(A^*),
\]
the restriction of $S^*(t)$ to $D(A^*)$ is a $C_0$-semigroup with
respect to the graph norm. Consequently, there exists a constant
$C_T>0$ such that
\[
\|\theta(0)\|_{D(A^*)}
=
\|S^*(\tau)\theta_\tau\|_{D(A^*)}
\leq
C_T\|\theta_\tau\|_{D(A^*)}.
\]
Therefore,
\begin{equation}\label{eq:tr_aux_1}
\left|
\left\langle
u_0,\theta(0)
\right\rangle_{D(A^*)',D(A^*)}
\right|
\leq
C_T
\|u_0\|_{D(A^*)'}
\|\theta_\tau\|_{D(A^*)}.
\end{equation}
Moreover, by the Cauchy--Schwarz inequality,
\[
\left|
\int_0^\tau\int_0^L
h(y,t)\theta_x(L,y,t)\,dy\,dt
\right|
\leq
\|h\|_{L^2((0,\tau)\times(0,L))}
\,
\|\theta_x(L,\cdot,\cdot)\|_{L^2((0,\tau)\times(0,L))}.
\]

To estimate the boundary trace, we rewrite the adjoint problem as a
forward evolution by reversing both the spatial and temporal variables. Then, we introduce the reflection operator $(R\varphi)(x,y)=\varphi(L-x,y),$
and notice that
\begin{equation*}
\partial_xR=-R\partial_x,
\qquad
\partial_x^{-1}R=-R(\partial_x^{-1})^*.
\end{equation*}
Indeed,  the first identity is the chain rule. For the second one, recalling that
$\partial_x^{-1}\varphi(x,y)=-\int_x^L\varphi(s,y)\,ds$ and
$(\partial_x^{-1})^{\ast}\varphi(x,y)=\int_0^x\varphi(s,y)\,ds$, the change of
variable $\sigma=L-s$ gives
\begin{equation*}
	\partial_x^{-1}(R\varphi)(x,y)
	=-\int_x^L\varphi(L-s,y)\,ds
	=-\int_0^{L-x}\varphi(\sigma,y)\,d\sigma
	=-\big[(\partial_x^{-1})^{\ast}\varphi\big](L-x,y),
\end{equation*}
that is, $\partial_x^{-1}R=-R(\partial_x^{-1})^{\ast}$. Composing this identity
with $R$ on the left and using $R^2=I$ we obtain also
$(\partial_x^{-1})^{\ast}R=-R\partial_x^{-1}$.

We then define
\begin{equation*}
w(x,y,s):=\theta(L-x,y,\tau-s),
\qquad (x,y,s)\in\Omega\times(0,\tau).
\end{equation*}
In particular, $w$ satisfies the corresponding forward system.
Moreover,
\begin{equation*}
w_x(0,y,s)=-\theta_x(L,y,\tau-s),
\end{equation*}
and therefore
\begin{equation*}
\|\theta_x(L,\cdot,\cdot)\|_{L^2((0,\tau)\times(0,L))}
=
\|w_x(0,\cdot,\cdot)\|_{L^2((0,\tau)\times(0,L))}.
\end{equation*}
Applying Proposition~\ref{pr:Kato} to $w$, precisely, from the estimate~\eqref{eq:traces} we obtain
\begin{equation*}
\|\theta_x(L,\cdot,\cdot)\|_{L^2((0,\tau)\times(0,L))}
\leq
C_T\|w(\cdot,\cdot,0)\|_{L^2(\Omega)}
=
C_T\|\theta_\tau\|_{L^2(\Omega)}
\leq
C_T\|\theta_\tau\|_{D(A^*)}.
\end{equation*}
Hence,
\begin{equation}\label{eq:tr_aux_2}
\left|
\int_0^\tau\int_0^L
h(y,t)\theta_x(L,y,t)\,dy\,dt
\right|
\leq
C_T
\|h\|_{L^2((0,T)\times(0,L))}
\|\theta_\tau\|_{D(A^*)}.
\end{equation}

Combining~\eqref{eq:tr_aux_1} and~\eqref{eq:tr_aux_2}, we conclude that
\[
|L_\tau(\theta_\tau)|
\leq
C_T
\left(
\|u_0\|_{D(A^*)'}
+
\|h\|_{L^2((0,T)\times(0,L))}
\right)
\|\theta_\tau\|_{D(A^*)}.
\]
Thus, $L_\tau$ is a continuous linear functional on $D(A^*)$. Then, by definition of the dual space, there exists a unique element $u(\tau)\in D(A^\ast)'$ represented by $L_\tau$. Hence, by Definition~\ref{transp-R},  the transposition identity
\eqref{transposition} holds. Taking the supremum over all
$\theta_\tau\in D(A^*)$ satisfying
$\|\theta_\tau\|_{D(A^*)}\leq1$, we obtain
\[
\|u(\tau)\|_{D(A^*)'}
\leq
C_T
\left(
\|u_0\|_{D(A^*)'}
+
\|h\|_{L^2((0,T)\times(0,L))}
\right).
\]
Since the constant is independent of $\tau\in[0,T]$, estimate
\eqref{eq:nh_est} follows.
\end{proof}

\section{Boundary observations}\label{Sec2}
\subsection{Observation for controllability}
It is well known in control theory \cite{Komornik,Lions88} that the exact controllability of a system is equivalent to proving an observability inequality.  Taking into account the following observability will ensure the controllability of the KP-II system \eqref{eq:KP}-\eqref{eq:BC}.
\begin{theorem}\label{th:Obs}
Let $L\in(0,+\infty)\setminus \mathcal{R}$ and $T>0$. Then there exists a constant $C:=C(T,L) > 0$ such that $\theta_T(x,y) \in L^2(\Omega)$ satisfies
\begin{equation}\label{eq:Obs}
\| \theta_T \|_{L^2(\Omega)}^2 \leq C \int_0^T\int_0^L \left|\theta_x(L,y,t) \right|^2 \,dy\,dt ,
\end{equation}
where $\theta$ is the solution of~\eqref{eq:tbk_h.1} with final data $\theta_T(x,y)$. 
\end{theorem}
\begin{proof} To prove this, first, making the change of variables $t\mapsto T-t$, $x\mapsto L-x$, transform~\eqref{eq:tbk_h.1} into the homogeneous problem~\eqref{eq:KP_h}. Then,~\eqref{eq:Obs} is equivalent to
\begin{equation}\label{eq:Obs.2}
\| u_0 \|_{L^2(\Omega)}^2 \leq C \int_0^T\int_0^L \left| u_x(0,y,t) \right|^2\,dy\,dt.
\end{equation}
To establish the observability inequality~\eqref{eq:Obs.2} we proceed in several steps.

\vspace{0.2cm}
\noindent \emph{\textbf{First step.}} By contradiction, assume that~\eqref{eq:Obs.2} does not hold. Then, there exists a sequence $\left(u_0^n\right)_n\subset L^2(\Omega)$ such that 
\begin{equation}\label{eq:Obs.3}
1 = \left\| u_0^n \right\|_{L^2(\Omega)}^2 > n  \int_0^T \int_0^L \left(u_x^n(0,y,t)\right)^2 \,dy\,dt,
\end{equation} 
for all $n\in \mathbb{N}$. Observe that \eqref{eq:Obs.3} is equivalent to the next two assertions: $\left\| u_0^n \right\|_{L^2(\Omega)}^2 = 1,$ and 
\begin{align}
&\int_0^T \int_0^L \left(u_x^n(0,y,t)\right)^2\,dy\,dt \longrightarrow 0,\quad\text{in } L^2((0,T)\times(0,L)) \label{eq:Obs.3-2},
\end{align}
where $u^n$ is the solution of~\eqref{eq:KP_h} with initial data $u_0^n$ for each $n\in\mathbb{N}$. From~\eqref{eq:reg} and~\eqref{eq:Obs.3} we get that $\left(u^n\right)_n\text{ is bounded in } L^2(0,T; H_x^1(\Omega)).$  Then, the sequence defined by 
$$u_t^n = -u_x^n - u_{xxx}^n - {\partial_x^{-1}}u_{yy}^n$$ is bounded in  $L^2(0,T; H^{-2}(\Omega)).$ Since $u^n(\cdot,\cdot,t) \in H_x^1(\Omega)$, thus $u_x^n(\cdot,\cdot,t)\in L^2(\Omega) \subset H^{-2}(\Omega)$. Moreover, thanks to the fact that $u_x^n(\cdot,\cdot,t)\in L^2(\Omega)$ we get that $u_{xxx}^n(\cdot,\cdot,t) \in H^{-2}(\Omega)$. 

\vspace{0.2cm}

\noindent \textbf{Claim 1.}  The sequence $\left(\partial_x^{-1}u_{yy}^n\right)_n$ is bounded in $L^2(0,T;H^{-2}(\Omega))$.

Indeed, since  $\partial_x^{-1}u_{yy}^n = \partial_y\left(\partial_x^{-1}u_y^n\right) $, it is enough to estimate this term in $H^{-1}(\Omega)$. Let $\xi\in H_0^1(\Omega)$. Then, by the definition of distributional
derivatives,
$$
\left\langle \partial_x^{-1}u_{yy}^n,\xi\right\rangle_{H^{-1}(\Omega),H_0^1(\Omega)}
= -\int_\Omega \left(\partial_x^{-1}u_y^n\right)\xi_y\,dxdy .
$$
Therefore, by the Cauchy--Schwarz inequality,
$$
\left\lvert\left\langle \partial_x^{-1}u_{yy}^n,\xi\right\rangle_{H^{-1}(\Omega),H_0^1(\Omega)}
\right\rvert
	\leq
\left\|\partial_x^{-1}u_y^n\right\|_{L^2(\Omega)} \|\xi_y\|_{L^2(\Omega)}.
$$
Since $\|\xi_y\|_{L^2(\Omega)}\leq\|\xi\|_{H_0^1(\Omega)},$ we obtain
$$
\left\|\partial_x^{-1}u_{yy}^n(\cdot,\cdot,t)\right\|_{H^{-1}(\Omega)}
	\leq
\left\|\partial_x^{-1}u_y^n(\cdot,\cdot,t)\right\|_{L^2(\Omega)}.
$$
In particular, due to the continuous embedding $H^{-1}(\Omega)\hookrightarrow H^{-2}(\Omega),$
there exists a constant $C>0$ such that 
$$
\left\|\partial_x^{-1}u_{yy}^n(\cdot,\cdot,t)\right\|_{H^{-2}(\Omega)}
	\leq
C\left\|\partial_x^{-1}u_y^n(\cdot,\cdot,t)\right\|_{L^2(\Omega)}.
$$
Integrating in time and exploiting the smoothing estimate obtained in~\eqref{eq:reg}, that is,
 $\partial_x^{-1}u_y^n\in L^2(0,T;L^2(\Omega)),$ we conclude that
$$
\left\|\partial_x^{-1}u_{yy}^n\right\|_{L^2(0,T;H^{-2}(\Omega))}^2
	\leq
C\int_0^T\left\|\partial_x^{-1}u_y^n(\cdot,\cdot,t)\right\|_{L^2(\Omega)}^2\,dt
	<\infty.
$$
Hence, $\partial_x^{-1}u_{yy}^n \in L^2(0,T;H^{-2}(\Omega)),$ which proves the claim.

It is important to emphasize that $H_x^1(\Omega)\hookrightarrow L^2(\Omega)$ is continuous but not compact. Therefore, the embeddings 
\begin{equation*}
H_x^1(\Omega)\hookrightarrow L^2(\Omega)\hookrightarrow H^{-2}(\Omega)
\end{equation*}
are not sufficient to directly apply the Aubin--Lions Lemma.

To tackle this problem, we exploit the smoothing effect provided by~\eqref{eq:reg}, namely,
\begin{equation}\label{eq:smoothing}
\int_0^T\!\int_{\Omega}\left(\partial_x^{-1}u_y^n\right)^2\,dx\,dy\,dt
\leq C\|u_0^n\|_{L^2(\Omega)}^2.
\end{equation}
Define the space
\begin{equation*}
\mathcal V
:=\left\{v\in L^2(\Omega)\; ;\; v_x\in L^2(\Omega)
\ \text{and }\ \partial_x^{-1}v_y\in L^2(\Omega)\right\},
\end{equation*}
endowed with the norm
\begin{equation*}
\|v\|_{\mathcal V}^2
:=\|v\|_{L^2(\Omega)}^2+\|v_x\|_{L^2(\Omega)}^2
+\|\partial_x^{-1}v_y\|_{L^2(\Omega)}^2.
\end{equation*}
Since $(u^n)_n$ is bounded in $L^2(0,T;H_x^1(\Omega))$ and satisfies the smoothing estimate
\eqref{eq:smoothing}, it follows that
\begin{equation*}
(u^n)_n \ \text{is bounded in } L^2(0,T;\mathcal V).
\end{equation*}

\medskip
\noindent\textbf{Claim 2.} The embedding $\mathcal V\hookrightarrow L^2(\Omega)$ is compact.
\medskip

{
Indeed, by the Fréchet--Kolmogorov theorem (see, for instance, \cite[Theorem~4.26]{Brezis0}), it suffices to show that translations are uniformly small on bounded subsets of $\mathcal V$.
For $(h,k)\in\mathbb R^2$, define
$$
\Omega_{h,k} := \left\{ (x,y)\in\Omega: (x+h,y+k)\in\Omega \right\},
$$
and
$$
\tau_{h,k}v(x,y) := v(x+h,y+k), \qquad (x,y)\in\Omega_{h,k}.
$$
Observe that, for $|h|,|k|<L$,
$$
\Omega_{h,k} = (0,L-|h|)\times(0,L-|k|).
$$
Hence,
$$
\Omega\setminus\Omega_{h,k} \subset \left( (L-|h|,L)\times(0,L)\right) \cup \left((0,L)\times(L-|k|,L)\right),
$$
and therefore $|\Omega\setminus\Omega_{h,k}|\leq L|h|+L|k|.$ 
In particular,
$$
|\Omega\setminus\Omega_{h,k}|\longrightarrow 0, \qquad\text{as }(h,k)\to(0,0).
$$

Now let $\mathfrak B\subset\mathcal V$ be bounded. Since $\mathcal{V}$ is continuously embedded in $L^2(\Omega)$, there exists $C>0$ such that
$$
\|v\|_{L^2(\Omega)}
	\leq
C\|v\|_{\mathcal V},
\qquad
\forall v\in\mathcal V.
$$
Consequently,
$$
\sup_{v\in\mathfrak B}
\|v\|_{L^2(\Omega)}
<\infty.
$$
Using the absolute continuity of the Lebesgue integral, we infer that
$$
\sup_{v\in\mathfrak B}
\int_{\Omega\setminus\Omega_{h,k}}
|v(x,y)|^2\,dx\,dy
\longrightarrow 0,
\qquad\text{as }(h,k)\to(0,0).
$$
Therefore, near the boundary it becomes zero for small translations.

Consequently, it is enough to estimate translations inside
$\Omega_{h,k}$.
We write
\begin{equation*}
\tau_{h,k}v-v
=
(\tau_{h,0}v-v)
+
\tau_{h,0}(\tau_{0,k}v-v),
\end{equation*}
so it is enough to estimate translations in the $x$ and $y$ directions separately.

Using the fundamental theorem of calculus, we have
\begin{equation*}
v(x+h,y)-v(x,y)
=
\int_0^h v_x(x+s,y)\,ds,
\end{equation*}
which implies
\begin{equation*}
\|\tau_{h,0}v-v\|_{L^2(\Omega_{h,0})}
\leq
|h|\,\|v_x\|_{L^2(\Omega)}.
\end{equation*}

Next, define $w:=\partial_x^{-1}v_y\in L^2(\Omega),$ so that $v_y=w_x$. Then,
\begin{equation*}
v(x,y+k)-v(x,y) = \int_0^k v_y(x,y+s)\,ds =
\partial_x\left(\int_0^k w(x,y+s)\,ds\right) =:\partial_x W_k(x,y).
\end{equation*}
For any $\varphi\in H_x^1(\Omega)$, integration by parts in $x$ yields
\begin{equation*}
\langle \tau_{0,k}v-v,\varphi\rangle_{H_x^{-1},H_x^1}
=
-\int_\Omega W_k\,\partial_x\varphi\,dx\,dy.
\end{equation*}
Hence,
\begin{equation*}
\|\tau_{0,k}v-v\|_{H_x^{-1}(\Omega)}
\leq
\|W_k\|_{L^2(\Omega)}
\leq
|k|\,\|\partial_x^{-1}v_y\|_{L^2(\Omega)}.
\end{equation*}

On the other hand,
\begin{equation*}
\|\tau_{0,k}v-v\|_{H_x^1(\Omega)}
\leq
2\|v\|_{H_x^1(\Omega)}.
\end{equation*}
Since $H_x^1(\Omega)\subset L^2(\Omega)\subset H_x^{-1}(\Omega)$ forms a Gelfand triple, interpolation yields the estimate
\begin{equation*}
\|f\|_{L^2(\Omega)}^2
\leq 
C\,\|f\|_{H_x^{-1}(\Omega)}\,\|f\|_{H_x^1(\Omega)}.
\end{equation*}
Applying this estimate, we infer that
\begin{equation*}
\|\tau_{0,k}v-v\|_{L^2(\Omega)}
\leq
C\,|k|^{1/2}\, \|v\|_{H_x^1(\Omega)}^{1/2} \|\partial_x^{-1}v_y\|_{L^2(\Omega)}^{1/2}.
\end{equation*}

Combining the previous estimates, we conclude that
\begin{equation*}
\|\tau_{h,k}v-v\|_{L^2(\Omega_{h,k})}
\longrightarrow 0,
\qquad\text{as }(h,k)\to(0,0),
\end{equation*}
uniformly for $v$ in bounded subsets of $\mathcal V$.} Therefore,
$\mathcal V$ is compactly embedded in $L^2(\Omega)$, which proves the claim 2.

\vspace{0.2cm}

{
Finally, since $\mathcal V\hookrightarrow L^2(\Omega)\hookrightarrow H^{-2}(\Omega),$ $(u^n)_n$ is bounded in $L^2(0,T;\mathcal V)$, and $(u_t^n)_n$ is bounded in $L^2(0,T;H^{-2}(\Omega))$, the Aubin--Lions theorem implies that $(u^n)_n$ is relatively compact in $L^2(0,T;L^2(\Omega))$.

Consequently, there exist a subsequence, still denoted by $(u^n)_n$, and a function $u\in L^2(0,T;L^2(\Omega))$ such that
$$
u^n\longrightarrow u, \qquad\text{strongly in }L^2(0,T;L^2(\Omega)).
$$

Moreover, since $(u^n)_n$ is bounded in $L^\infty(0,T;L^2(\Omega))$ and $(u_t^n)_n$ is bounded in
$L^2(0,T;H^{-2}(\Omega))$, a standard compactness result (see~\cite[Corollary~4]{SimonCompact}) implies, after extracting a subsequence, that
\begin{equation*}
u^n\longrightarrow u \quad\text{in}\quad C([0,T];H^{-2}(\Omega)).
\end{equation*}
In particular,
\begin{equation*}
u^n(\cdot,\cdot,t) \longrightarrow u(\cdot,\cdot,t)
\quad\text{in}\quad H^{-2}(\Omega),
\qquad \forall\,t\in[0,T].
\end{equation*}
}

\vspace{0.2cm}
\noindent \textbf{Claim 3.} { The nonlocal term satisfies
$$
\partial_x^{-1}u_{yy}^n \longrightarrow \partial_x^{-1}u_{yy}
\quad\text{in}\quad L^2(0,T;H^{-2}(\Omega)).
$$}
In fact, from the definition of the operator $\partial_x^{-1}$ we have that ${\partial_x^{-1}} u^n = \varphi^n$ where $\varphi_x^n = u^n$, $u^n(\cdot,\cdot,t) \in H_x^1(\Omega)$ and $\varphi^n(\cdot,\cdot,t) \in H_x^1(\Omega)$. Since ${\partial_x^{-1}}u_{yy}^n = \varphi_{yy}^n$ it follows that
\begin{equation*}
\begin{aligned}
\left\| {\partial_x^{-1}} u_{yy}^n(\cdot,\cdot,t) - {\partial_x^{-1}} u_{yy}(\cdot,\cdot,t)  \right\|_{H^{-2}(\Omega)}
	& \leq C \| \varphi^n(\cdot,\cdot,t) - \varphi(\cdot,\cdot,t) \|_{L^{2}(\Omega)}\\
	& \leq C L^2 \| \varphi_{x}^n(\cdot,\cdot,t) - \varphi_{x}(\cdot,\cdot,t) \|_{L^{2}(\Omega)}\\
	& = C L^2 \| u^n(\cdot,\cdot,t) - u(\cdot,\cdot,t) \|_{L^{2}(\Omega)} \longrightarrow 0,
\end{aligned}
\end{equation*}
showing claim 3. 
\vspace{0.1cm}

It follows from~\eqref{eq:Time_est} and~\eqref{eq:Obs.3} that $(u_0^n)_n$ is a Cauchy sequence in $L^2(\Omega)$, hence  it is {\ convergent in} $L^2(\Omega)$. Let $u_0:=\lim_{n\to\infty}u_0^n$. Then, $u^n=S(\cdot)u_0^n\to S(\cdot)u_0$ in $L^2(0,T;L^2(\Omega))$ (so that $u=S(\cdot)u_0$), with
\begin{equation}\label{eq:Obs.3-1}
\| u_0 \|_{L^2(\Omega)} = 1.
\end{equation}
From~\eqref{eq:Obs.3} follows that \color{black}
\begin{equation*}
0  = \liminf_{n\to+\infty} \left[ \int_0^T \int_0^L \left(u_x^n(0,y,t)\right)^2 \,dy\,dt\right] \geq \int_0^T \int_0^L \left(u_x(0,y,t)\right)^2 \,dy\,dt.
\end{equation*}
Then, $ u_x(0, \cdot ,\cdot)  = 0$ a.e. in $(0,L)\times (0,T)$.  Hence, taking into account the previous claims, $u$ is a solution for
\begin{equation*}
\begin{cases}
	u_t + u_x + u_{xxx} + \partial_x^{-1} u_{yy}= 0,\\
	u(0,y) = u(L,y) = u(x,0) = u(x, L) = u_{x}(L,y) = u_{x}(0,y) = 0,\\
	u(x,y,0) = u_0(x,y),
\end{cases}
\end{equation*}
satisfying \eqref{eq:Obs.3-1}.

\vspace{0.2cm}
\noindent \emph{\textbf{Second step.}} We reduce the proof of the observability inequality into a spectral problem as done in \cite{Rosier} for the KdV equation. The result is the following.

\begin{lemma}\label{LemmaRosier}
For any $L\in(0,+\infty)\setminus \mathcal{R} $ and $T > 0$, let $N_T$ denote the space of all the initial states,  $u_0 \in L^2(\Omega)$ for which the solution $u(t)=S(t)u_0$ of \eqref{eq:KP_h} satisfies $u_x(0,y) = 0$. Then  $N_T =\{0\}$. 
\end{lemma}
\begin{proof}
Using the arguments as those given in \cite[Lemma 3.4]{Rosier}, it follows that if $N_T \neq \{0\}$, the map $u_0 \in N_T \mapsto A(N_T) \subset \mathbb{C}N_T$ (where $\mathbb{C}N_T$ denote the complexification of $N_T$) has (at least) one eigenvalue; hence, there exist $\lambda \in \mathbb{C}$ and  $u_0 \in H^3_x(\Omega)\cap H^2_y(\Omega)$\footnote{where $H_y^k(\Omega):=
	\left\lbrace%
	\varphi \colon \partial_y^j \varphi \in L^2(\Omega),\ \text{for } 0\leq j \leq k
	\right\rbrace ,$
	equipped with the norm $\left\lVert \varphi \right\rVert_{H_y^k(\Omega)}^2 =\sum_{j=0}^k \left\lVert \partial_y^j \varphi\right\rVert_{L^2(\Omega)}^2.$ } such that 
\begin{equation}\label{eq:KP_sp}
\begin{cases}
		\lambda u_0(x,y) + u_{x,0}(x,y) + u_{xxx,0}(x,y) + \partial_x^{-1} u_{yy,0}(x,y)= 0,\\
	u_0(0,y) = u_0(L,y) = u_0(x,0) = u_0(x, L) = u_{x,0}(L,y)= u_{x,0}(0,y) = 0,
\end{cases}
\end{equation}
with $(x,y) \in \Omega$. To conclude the proof of the lemma, we will prove that this does not hold.
\end{proof}

\noindent \emph{\textbf{Third step.}} To obtain the contradiction, it remains to prove that a pair $(\lambda, u_0)$ as above does not exist.  Precisely, this step is to show that no nontrivial solution exists for the spectral problem \eqref{eq:KP_sp}.
\begin{lemma}\label{le:RedSPaa}
Let $L\in(0,+\infty)$. Consider the following assertion:
$$
\text{$(\mathcal{F}) \quad \exists \lambda \in \mathbb{C}, \exists u_0 \in H^3_x(\Omega) \cap H^2_y(\Omega) \backslash\{0\}$ fulfilling \eqref{eq:KP_sp}.}
$$
Then, $(\mathcal{F})$ holds if and only if $L \in \mathcal{R}$.
\end{lemma}
\begin{proof}
To simplify the notation, let us denote $u(x,y) = u_0(x,y)$. Then~\eqref{eq:KP_sp} can be rewritten as 
\begin{equation*}
\begin{cases}
\lambda u + u_x + u_{xxx} + {\partial_x^{-1}}u_{yy} = 0,\\
u(0,y) = u(L,y) = u_x(0,y) = u_x(L,y) =u(x,0) = u(x,L) = 0.
\end{cases}
\end{equation*}
\color{black} 


First, let us focus on the structure in the $y$-direction. Let $\{\varphi_n\}_{n\geq 1}$ be the orthonormal basis of $L^2(0,L)$ given by
\begin{equation*}
-\varphi_n''(y)=\mu_n\varphi_n(y),\qquad \varphi_n(0)=\varphi_n(L)=0,
\qquad 
\mu_n=\left(\frac{n\pi}{L}\right)^2,
\end{equation*}
that is,
\begin{equation*}
\varphi_n(y)=\sqrt{\frac{2}{L}}\sin\left(\frac{n\pi y}{L}\right),\qquad \forall n\ge 1.
\end{equation*}
Expanding $u(x,y)$ as
\begin{equation*}
u(x,y)=\sum_{n\geq 1} u_n(x)\,\varphi_n(y),
\qquad \text{with}\ \ 
u_n(x):=\int_0^L u(x,y)\,\varphi_n(y)\,dy,
\end{equation*}
we get
\begin{equation*}
u_{yy}(x,y)=\sum_{n\geq 1} u_n(x)\,\varphi_n''(y)
=-\sum_{n\geq 1}\mu_n\,u_n(x)\,\varphi_n(y).
\end{equation*}
Since $\partial_x^{-1}$ acts only on the $x$--variable, it commutes with the projection in $y$, and therefore
\begin{equation*}
\partial_x^{-1}u_{yy}(x,y)
=-\sum_{n\geq 1}\mu_n\,\partial_x^{-1}u_n(x)\,\varphi_n(y).
\end{equation*}
Substituting these representations into the eigenvalue equation
\begin{equation*}
\lambda u + u_x + u_{xxx} + \partial_x^{-1}u_{yy}=0,
\end{equation*}
and using the orthonormality of $\{\varphi_n\}_{n\geq 1}$, we infer that each coefficient $u_n$ satisfies
\begin{equation}\label{eq:mode_eq_y}
\lambda u_n(x)+u_n'(x)+u_n'''(x)-\mu_n\,\partial_x^{-1}u_n(x)=0,
\qquad \forall n\geq 1.
\end{equation}
%

To simplify the notation, we denote $u_n(x)$ simply by $p(x)$ for all $n\geq 1$. Now, assume that there exists $p\in H^3(0,L)\setminus \lbrace 0 \rbrace$ solution of~\eqref{eq:mode_eq_y}, thus, 
$p$ satisfies 
\begin{equation}\label{eq:O.KdV.1}
\begin{cases}
\lambda p(x) + p'(x) + p'''(x) - \left(\frac{n\pi}{L}\right)^2 {\partial_x^{-1} p(x)} = 0 \\
p(0) = p(L) = p'(0) = p'(L) = 0.
\end{cases}
\end{equation}\color{black}
Denote by $$\hat{p}(k) = \int_0^L e^{-ixk} p(x)dx.$$ Multiplying~\eqref{eq:O.KdV.1} by $ike^{-ixk}$, integrating by parts over $(0,L)$ and using the boundary conditions we obtain 
\begin{equation*}
\begin{aligned}
-\left(\frac{n\pi}{L}\right)^2 \hat{p}(k) + \lambda ik\hat{p}(k) + (ik)^2\hat{p}(k) + (ik)^4\hat{p}(k) \\
	\qquad = \left(\frac{n\pi}{L}\right)^2 {\partial_x^{-1}}p(0) + (ik)p''(0) - (ik)e^{-iLk} p''(L).
\end{aligned}
\end{equation*}
Indeed, note first that
\begin{equation*}
\begin{aligned}
(ik)\widehat{\partial_x^{-1}p}(k) 
	& =(ik) \int_0^L e^{-ixk} \partial_x^{-1}p(x)\,dx \\
	 & = - \int_0^L \left[ e^{-ixk}\right]_x {\partial_x^{-1}}p(x)\,dx 
	 = \hat{p}(k) +  {\partial_x^{-1}}p(0)
\end{aligned}
\end{equation*}
and
\begin{equation*}
\begin{aligned}
(ik)\widehat{p'''}(k) 
	& = (ik)\int_0^L e^{-ixk} p'''(x)\,dx 
	= (ik)^4 \hat{p}(k) + (ik)e^{-iLk}p''(L) - (ik)p''(0).
\end{aligned}
\end{equation*}
Due to the boundary conditions, $p$ and $p'$  have traces equal to zero. Consequently,
\begin{equation*}
\left(\lambda(ik) + (ik)^2 + (ik)^4 - \left(\frac{n\pi}{L}\right)^2 \right) \hat{p}(k) = \left(\frac{n\pi}{L}\right)^2 {\partial_x^{-1}}p(0) + (ik) p''(0) - (ik)e^{-iLk} p''(L).
\end{equation*}
By rearranging the equation above setting $\lambda = i\sigma \in \mathbb{C}$, follows that
\begin{equation*}
\hat{p}(k) = \frac{ \left(\frac{n\pi}{L}\right)^2 {\partial_x^{-1}}p(0) + (ik) p''(0) - (ik) e^{-iLk} p''(L) }{k^4 - k^2 - \sigma k - \left(\frac{n\pi}{L}\right)^2}.
\end{equation*}

Using the Paley–Wiener theorem (see \cite[Section 4, p. 161]{Yosida}) and the usual characterization of $H^3(\mathbb{R})$ by means of the Fourier transform we see that $\mathcal{F}$ is equivalent to the existence of $\sigma$ and $\alpha_0,\alpha_1,\alpha_2$ such that the function
\begin{equation}\label{eq:R0}
f(k) := \frac{\left(\frac{n\pi}{L}\right)^2\alpha_0 + (ik) \alpha_1 - (ik) e^{-iLk} \alpha_2}{k^4 - k^2 - \sigma k - \left(\frac{n\pi}{L}\right)^2} = \frac{R(k)}{Q(k)}
\end{equation}
satisfies
\begin{enumerate}
\item[(i)] $f$ is an entire function in $\mathbb{C}$;
\item[(ii)] $\displaystyle\int_{\mathbb{R}} |f(k)|^2 \left(1+|k|^2\right)^2\,dk < \infty$;
\item[(iii)] There exist positive constants $C, N$ such that for all $k\in\mathbb{C}$, $$|f(k)| \leq C(1+|k|)^N e^{L |\text{Im}\, k|}.$$
\end{enumerate}

 Let $\kappa_0,\kappa_1,\kappa_2,\kappa_3$ be the roots of $Q(k)$ and consider two cases:
\begin{itemize}
     \item {\bf Case 1}: Suppose that  $\alpha_0= 0$. 
     \end{itemize} 
     If it is the case, the related trace value ${\partial_x^{-1}} p(0) = 0$ and~\eqref{eq:R0} turns into
\begin{equation}\label{eq:R5}
f(k) = \frac{ (ik) \alpha_1 - (ik) e^{-iLk} \alpha_2}{k^4 - k^2 - \sigma k - \left(\frac{n\pi}{L}\right)^2}= \frac{R_1(k)}{Q(k)}.
\end{equation}
Note that $Q(0) = - \left(\frac{n\pi}{L}\right)^2 \neq0$ and consequently $k=0$ cannot be a root of $R_1(k)$ and $Q(k)$ simultaneously. Since the roots of $R_1(k)$ are simple, unless $\alpha_1 = \alpha_2 = 0$, (i) holds, then (ii) and (iii) are satisfied. It follows that our problem is equivalent to the existence of complex numbers $\sigma, \kappa_0 \in \mathbb{C}$ and  $m_1, m_2, m_3\in\mathbb{N}^*$
such that if we set
$$
\kappa_1 := \kappa_0 + m_1\frac{2\pi}{L},\quad \kappa_2 := \kappa_1 + m_2\frac{2\pi}{L}\quad\text{and}\quad \kappa_3 := \kappa_2 + m_3\frac{2\pi}{L}
$$
we have
$$
Q(k) = (k - \kappa_0)(k - \kappa_1)(k - \kappa_2)(k - \kappa_3).
$$
By using the Girard–Newton formula, we get
\begin{equation}\label{eq:R0.1}
\begin{cases}
0& =\sum_{j=0}^3 \kappa_j, \\
-1 &= \kappa_0\sum_{j=1}^3 \kappa_j  + \kappa_1 \sum_{j=2}^3 \kappa_j + \kappa_2\kappa_3 , \\
\sigma &=\kappa_0\kappa_1\sum_{j=2}^3 \kappa_j + \kappa_2\kappa_3\sum_{j=0}^1 \kappa_j, \\
 -\left(\frac{n\pi}{L}\right)^2  &=\prod_{j=0}^3 \kappa_j.
 \end{cases}
\end{equation}
Observe that, $\kappa_2, \kappa_3$ can be rewritten in terms of $\kappa_0$ as
\begin{equation*}
\kappa_2 = \kappa_0 + (m_1+m_2) \frac{2\pi}{L}\quad\text{and}\quad\kappa_3 = \kappa_0 + (m_1+m_2 + m_3) \frac{2\pi}{L}. 
\end{equation*}
From the first relation in \eqref{eq:R0.1}, it follows that
\begin{equation*}
\kappa_0 = \frac{-\pi (3m_1 + 2m_2 + m_3)}{2L}.
\end{equation*}

Now, let us introduce the notation 
\begin{equation}\label{eq:SH.1}
q=\frac{\pi}{2L},
\qquad
\mathfrak a = m_1+2m_2+m_3.
\end{equation}
Using the aforementioned notation, the roots can be written as
\begin{align}\label{eq:SH.2}
\begin{aligned}
	\kappa_0&=-q(\mathfrak{a}+2m_1),&
	\kappa_1&=q(2m_1-\mathfrak{a}),\\
	\kappa_2&=q(\mathfrak{a}-2m_3),&
	\kappa_3&=q(\mathfrak{a}+2m_3).
\end{aligned}
\end{align}

Therefore, employing~\eqref{eq:SH.2} in the second Girard-Newton relation~\eqref{eq:R0.1} and performing straightforward computations, we get
\begin{align*}
	-1
	={}&-q^2(\mathfrak a+2m_1)^2
	+2q^2\mathfrak a(2m_1-\mathfrak a)
	+q^2(\mathfrak a^2-4m_3^2)\\
	={}&q^2\left[
	-(\mathfrak a+2m_1)^2
	+2\mathfrak a(2m_1-\mathfrak a)
	+\mathfrak a^2-4m_3^2
	\right]\\
	={}& -2q^2\left(\mathfrak a^2+2m_1^2+2m_3^2\right),
\end{align*}
Since $q=\pi/2L$, it follows that
\begin{equation}\label{eq:SH.3}
	L
	=
	\frac{\pi}{\sqrt{2}}
	\sqrt{
		\mathfrak{a}^2+2m_1^2+2m_3^2
	}.
\end{equation}

On the other hand, the fourth relation in~\eqref{eq:R0.1} gives
\begin{equation*}
	-\left(\frac{n\pi}{L}\right)^2
	=\prod_{j=0}^{3}\kappa_j
	=q^4
	\left(\mathfrak{a}^2-4m_1^2\right)
	\left(\mathfrak{a}^2-4m_3^2\right).
\end{equation*}
Combining this identity with~\eqref{eq:SH.3}, we obtain the compatibility condition
\begin{equation}\label{eq:compatibility-condition}
	n^2
	=-\frac{
	\left(
		\mathfrak{a}^2-4m_1^2
	\right)
	\left(
		\mathfrak{a}^2-4m_3^2
	\right)}{
	8\left(
		\mathfrak{a}^2+2m_1^2+2m_3^2
	\right)},
	\quad\text{and}\quad
	\mathfrak{a}=m_1+2m_2+m_3.
\end{equation}

Notice that the denominator in~\eqref{eq:compatibility-condition} is strictly positive. Moreover, the right-hand side is positive if and only if $|m_1-m_3|>2m_2$.
Therefore~\eqref{eq:compatibility-condition} hold for some $n\in\mathbb{N}^*$, whenever $(m_1,m_2,m_3)\in\mathcal{M}$, where $\mathcal{M}$ is defined in~\eqref{eq:cond-M}.

\begin{itemize}
    \item {\bf Case 2}: Suppose that  $\alpha_0\neq 0$.  
    \end{itemize}
For this case notice that $$R(k) = \frac{1}{2}\left(R(k) +  \overline{R(\bar{k})}\right) + \frac{1}{2}\left(R(k) - \overline{R(\bar{k})}\right),$$ where
$$\overline{R(\bar{k})} = \left(\frac{n\pi}{L}\right)^2\alpha_0 - (ik) \alpha_1 + (ik) e^{iLk} \alpha_2. $$
Therefore,
\begin{equation}\label{eq:R1}
R(k) + \overline{R(\bar{k})} 
	 = 2\left(\dfrac{n\pi}{L}\right)^2\alpha_0  + (ik)\alpha_2\left(e^{iLk} - e^{-iLk}\right)
	 = 2\left(\dfrac{n\pi}{L}\right)^2\alpha_0  - 2k\alpha_2\sin(Lk)
\end{equation}
and
\begin{equation}\label{eq:R2}
R(k) - \overline{R(\bar{k})} 
	 = 2(ik)\alpha_1  - (ik)\alpha_2\left(e^{iLk} + e^{-iLk}\right)
	 = 2 (ik)\alpha_1  - 2(ik)\alpha_2\cos(Lk).
\end{equation}
Here, we used that 
$$\cos(z) = \frac{e^{iz} + e^{-iz}}{2}\quad\text{and}\quad \sin(z) = \frac{e^{iz} - e^{-iz}}{2i}.$$
Rearranging~\eqref{eq:R1} and~\eqref{eq:R2} and gathering the expressions, we get
\begin{equation*}
 R(k) \Longleftrightarrow\quad\begin{cases}
\left(\dfrac{n\pi}{L}\right)^2\alpha_0  - k\alpha_2\sin(Lk) = 0, \\
\alpha_1 k - k\alpha_2\cos(Lk) = 0,
\end{cases}
\Longleftrightarrow \quad
\begin{cases}
   	k\alpha_2\sin(Lk) = \left(\dfrac{n\pi}{L}\right)^2\alpha_0, \\
 	k\alpha_2\cos(Lk) = \alpha_1 k,
\end{cases}
\end{equation*}
or, equivalently
\begin{equation}\label{eq:R3}
\begin{cases}
   	k^2\alpha_2^2\sin^2(Lk) = \left(\dfrac{n\pi}{L}\right)^4\alpha_0^2,\\
 	k^2\alpha_2^2\cos^2(Lk) = \alpha_1^2 k^2 .
\end{cases}
\end{equation}
Adding the results in~\eqref{eq:R3}, follows that
\begin{equation}\label{eq:R3.1}
k^2(\alpha_2^2 - \alpha_1^2) = \alpha_0^2 \left(\dfrac{n\pi}{L}\right)^4\neq 0 .
\end{equation}
From \eqref{eq:R3.1}, it follows that 
$$
    \kappa_j^2 =\frac{\alpha_0^2}{(\alpha_2^2 - \alpha_1^2)}  \left(\dfrac{n\pi}{L}\right)^4, \quad \forall j =0,1,2,3.
$$
Thus, if $\alpha^2_2-\alpha_1^2>0$, we have that $$\kappa_j =\pm \sqrt{\frac{\alpha_0^2}{(\alpha_2^2 - \alpha_1^2)}}  \left(\dfrac{n\pi}{L}\right)^2,$$ and if  $\alpha^2_2-\alpha_1^2<0$, it follows that $$\kappa_j =\pm i\sqrt{\frac{\alpha_0^2}{(\alpha_2^2 - \alpha_1^2)}}  \left(\dfrac{n\pi}{L}\right)^2.$$ In any case, without loss of generality, we deduce that
$$
Q(k)=(k-\kappa_0)^2(k-\kappa_2)^2.
$$

Note that $Q(k)$ (see \eqref{eq:R5}) cannot have two roots of order two. Indeed, note that
\begin{equation*}
\begin{aligned}
Q(k) 
    & = k^4 - k^2 -\sigma k - \left(\dfrac{n\pi}{L}\right)^2 \\
    & = (k-\kappa_0)^2(k-\kappa_2)^2 \\
    & = k^4 -2(\kappa_0+\kappa_2) k^3 + (\kappa_2^2 + 4\kappa_0\kappa_2 + \kappa_0^2) k^2 -2(\kappa_0\kappa_2^2 + \kappa_0^2\kappa_2) k +\kappa_0^2\kappa_2^2.
\end{aligned}
\end{equation*}
Therefore, $\kappa_0$ and $\kappa_2$ must satisfy
\begin{equation*}
\begin{cases}
\kappa_0 + \kappa_2 = 0; \\
\kappa_0^2 + 4\kappa_0\kappa_2 + \kappa_2^2 = -1; \\
\kappa_0^2\kappa_2^2 = -\left(\dfrac{n\pi}{L}\right)^2.
\end{cases}
\end{equation*}
The first relation implies that $\kappa_0 = -\kappa_2$, then from the second relation we get $\kappa_0^2 = \frac{1}{2}.$ With this in hand, from the third relation, we obtain that $\frac{1}{4} = -\left(\frac{n\pi}{L}\right)^2$, giving a contradiction. Hence, $(\mathcal{F})$ holds if and only if $L\in\mathcal{R}$. This completes the proof of Lemma \ref{le:RedSPaa} and, consequently, the proof of Lemma \ref{LemmaRosier}.
\end{proof}
With these lemmas in hand, the proof of Theorem \ref{th:Obs} is achieved.
\end{proof}

\subsection{Observation for stabilization}  Due to the structure of the energy for the KP-II with one boundary feedback,  we have a slightly different observability inequality that will ensure the stabilization of the KP-II system \eqref{eq:KP}-\eqref{eq:BC.f}. 
\begin{theorem}
Let $L\in(0,+\infty)\setminus \mathcal{R}$ and $T>0$. Then there exists a constant $C:= C(T,L) > 0$ such that  $u_0(x,y) \in L^2(\Omega)$ satisfies
\begin{equation}\label{eq:S.Obs}
\|u_0\|_{L^2(\Omega)}^2
	\leq C	\left(	\int_0^T\int_0^L u_{x}^2(0,y,t)\,dy\,dt + \int_0^T\int_0^L({\partial_x^{-1}} u_y(0,y,t))^2\,dy\,dt 	\right),
\end{equation}
where $u$ is the solution of~\eqref{eq:KP}-\eqref{eq:BC.f} with initial data $u_0(x,y)$. 
\end{theorem}

\begin{proof}
Essentially, the proof of this result is similar to that done in the previous subsection. It is important to emphasize that in the proof of~\eqref{eq:S.Obs} we get an additional condition to reach the contradiction; however, the reasoning is the same. Recalling the spectral problem
\begin{equation*}
\begin{cases}
\lambda u + u_x + u_{xxx} + {\partial_x^{-1}}u_{yy} = 0\\
u(0,y) = u(L,y) = u_x(0,y) = u_x(L,y) = u(x,0) = u(x,L) = {\partial_x^{-1}}u(0,y) = 0.
\end{cases}
\end{equation*}
{\color{black}
we expand $u$ in the orthonormal basis $ \{\varphi_n\}_{n\geq 1}$ of $L^2(0,L)$ associated with the Dirichlet Laplacian in the $y$-variable
\begin{equation*}
u(x,y) = \sum_{n\geq 1} u_n(x)\varphi_n(y).
\end{equation*}

As shown in the proof of Lemma~\ref{le:RedSPaa}, this decomposition leads, for each $n\geq 1$, to the one-dimensional equation~\eqref{eq:mode_eq_y}. Once more, to simplify the notation, we denote $u_n(x)$ simply by $ p(x)$. Consequently, for each $n\geq 1$, the function $p$ satisfies 
\begin{equation}\label{eq:S.O.KdV.1}
\begin{cases}
\lambda p(x) + p'(x) + p'''(x) - \left(\dfrac{n\pi}{L}\right)^2 \partial_x^{-1} p(x) = 0,\\
p(0) = p(L) = p'(0) = p'(L) = \partial_x^{-1}p(0) = 0.
\end{cases}
\end{equation}

}

Now, by the same reasoning as in the proof of the Theorem~\ref{th:Obs}, denoting by $\hat{p}(k) = \int_0^L e^{-ixk} p(x)\,dx$, multiplying~\eqref{eq:S.O.KdV.1} by $ike^{-ixk}$, integrating by parts over $(0,L)$ and using the boundary conditions we obtain
$$
-\left(\frac{n\pi}{L}\right)^2  \hat{p}(k) + \lambda (ik) \hat{p}(k) + (ik)^2\hat{p}(k) + (ik)^4\hat{p}(k) 
	 = (ik)p''(0) - (ik)e^{-iLk} p''(L).
$$
Consequently, after  rearranging and setting $\lambda = i\sigma$,
\begin{equation*}
\hat{p}(k) =(ik) \frac{   p''(0) -  e^{-iLk} p''(L) }{k^4 - k^2 - \sigma k - \left(\frac{n\pi}{L}\right)^2}.
\end{equation*}
Using the Paley-Wiener theorem, we can characterize the Fourier transform by showing the existence of $\sigma$ and $\alpha_1,\alpha_2$ such that the application
\begin{equation}\label{eq:S.f}
f(k) := (ik)\frac{(  \alpha_1 -  e^{-iLk} \alpha_2)}{k^4 - k^2 - \sigma k - \left(\frac{n\pi}{L}\right)^2}
\end{equation}
satisfies the same properties of the function $f$ defined in the proof of Lemma \ref{le:RedSPaa}. Observe that~\eqref{eq:S.f} is equal to \eqref{eq:R5} and therefore ~\eqref{eq:S.Obs} holds.
\end{proof}

\section{Control results of KP-II equation}\label{Sec3}
\subsection{Boundary controllability}
To prove the control result, let us follow the H.U.M. developed by J.-L. Lions~\cite{Lions88}. Next, we define the exact controllability property.
\begin{definition}\label{df:ex_control}
Let $T>0$. We say that the system \eqref{eq:KP}-\eqref{eq:BC} is exactly controllable in time $T$ if for any initial and final data $u_0, u_T \in L^2(\Omega)$, there exists $h\in L^2\left((0,T)\times (0,L)\right)$ such that $u(\cdot,\cdot, T) = u_T(\cdot,\cdot)$.
\end{definition}

Now, we use the adjoint system to give an equivalent condition for Definition~\ref{df:ex_control}. The condition is the following one.
\begin{lemma}
Let  $u_0, u_T \in L^2(\Omega)$. Then there exists $h\in L^2\left((0,T)\times (0,L)\right) $, such that the solution $u$ of \eqref{eq:KP}-\eqref{eq:BC} satisfies $u(x,y,T) = u_T(x,y)$ if and only if
\begin{equation}\label{eq:Le_Co}
\int_{\Omega} \theta_T(x,y) u(x,y,T) \,dx\,dy =  \int_0^T \int_0^L h(y,t)\theta_x(L,y,t)\,dy\,dt + \int_{\Omega} \theta(x,y,0) u_0(x,y) \,dx\,dy,
\end{equation}
for any $\theta_T(x,y) \in L^2(\Omega)$ and $\theta$ being the solution of the adjoint problem~\eqref{eq:tbk_h.1}.
\end{lemma}

\begin{proof}[Proof of Theorem~\ref{th:Control}]
Notice that the relation \eqref{eq:Le_Co}  may be seen as an optimality condition for the critical points of the functional $\Lambda: L^2(\Omega)\to \mathbb{R}$, defined by
\begin{equation}\label{eq:J.op}
\Lambda(\theta_T) = \frac{1}{2}\left\| \theta_x(L,\cdot,\cdot)\right\|_{L^2((0,T)\times (0,L))}^2 - \int_{\Omega} \theta_T u(x,y,T)\,dx\,dy,
\end{equation}
where $\theta$ is the solution of~\eqref{eq:tbk_h.1}.  A control may be obtained from the solution of the homogeneous system \eqref{eq:tbk_h.1} with the initial data minimizing the functional $\Lambda$. Hence, the controllability is reduced to a minimization problem. To guarantee that $\Lambda$ defined by~\eqref{eq:J.op} has a unique minimizer, we use the next fundamental result in the calculus of variations.
\begin{theorem}[See~\cite{Brezis0}]
Let $H$ be a reflexive Banach space, $K$ a closed convex subset of $H$ and $\Lambda: K \rightarrow \mathbb{R}$ a function with the following properties:
\begin{enumerate}
\item $\Lambda$ is convex and lower semi-continuous;
\item If $K$ is unbounded then $\Lambda$ is coercive.  
\end{enumerate}
Then $\Lambda$ attains its minimum in $K$, i. e. there exists $x_0 \in K$ such that
$
\Lambda\left(x_0\right)=\min _{x \in K} \varphi(x).
$
\end{theorem}

Note that $\Lambda$ defined in~\eqref{eq:J.op} is continuous and convex. The existence of a minimum is ensured if we prove that $\Lambda$ is also coercive, which is obtained with the \emph{observability inequality} given by \eqref{eq:Obs}. Thus, 
$$
\begin{aligned}
\Lambda (\theta_T)  
	& = \frac{1}{2}\| \theta_x(L,\cdot,\cdot) \|_{L^2((0,T)\times(0,L))}^2 - \int_{\Omega} \theta_T u(x,y,T)\,dx\,dy \geq\frac{C^{-1}}{2} \|\theta_T\|_{L^2(\Omega)}^2.
\end{aligned}
$$
Therefore, the KP-II equation~\eqref{eq:KP}-\eqref{eq:BC} is exactly controllable.
\end{proof}

\subsection{Boundary stabilization}
We aim to show that the KP-II equation \eqref{eq:KP}-\eqref{eq:BC} can be stabilized by selecting an appropriate feedback-damping control law. Specifically, by choosing \( h(y,t) = -\alpha u_x(0,y,t) \) with \( 0 < |\alpha| < 1 \), we can establish our second main result. Recall that the energy of the KP-II equation, defined by \eqref{eq:S.En}, is a nonincreasing function due to \eqref{eq:S.En'}. 
\begin{proof}[Proof of the Theorem~\ref{th:S.Exp}]
Through the utilization of energy dissipation, specifically~\eqref{eq:S.E'(t)}, in conjunction with the \emph{observability inequality}~\eqref{eq:S.Obs}, we obtain the exponential stabilization. Indeed, notice that, integrating ~\eqref{eq:S.E'(t)} over $[0,T]$, it follows that
\begin{equation*}
E(T) - E(0)  \leq -C	\left[\int_0^T\int_0^L u_{x}^2(0,y,t)\,dy\,dt + \int_0^T\int_0^L({\partial_x^{-1}} u_y(0,y,t))^2\,dy \,dt \right].
\end{equation*}

Thus, we have that $$E(T) - E(0)  \leq -\frac{1}{C} E(0).$$ Since the energy is dissipative, it follows that $E(T)\leq E(0)$, thus $$E(T) \leq E(0) -\frac{1}{C} E(0) \leq E(0) - \frac{1}{C} E(T)\implies E(T) \leq \delta E(0), \ \text{where}\ \delta=\frac{C}{1+C}<1.$$
Now, applying the same argument on the interval $[(m-1) T, m T]$ for $m=1,2, \ldots$, yields that
$$
E(m T) \leq \delta E((m-1) T) \leq \cdots \leq \delta^m E(0).
$$
Thus, we have
$$
E(m T) \leq e^{-\mu_0 m T} E(0) \quad \text { with } \quad \mu_0=\frac{1}{T} \ln \left(1+\frac{1}{C}\right)>0 .
$$
For an arbitrary $t>0$, there exists $m \in \mathbb{N}^*$ such that $(m-1) T<t \leq m T$, and by the nonincreasing property of the energy, we conclude that
$$
E(t) \leq E((m-1) T) \leq e^{-\mu_0 (m-1) T} E(0) \leq \frac{1}{\delta} e^{-\mu_0 t} E(0),
$$
showing uniform exponential stability.
\end{proof}

\subsection{An explicit decay rate}
Finally, by using Lyapunov's approach, we can give an explicit decay rate for the solutions of \eqref{eq:KP}-\eqref{eq:BC.f}.
\begin{proof}[Proof of Theorem~\ref{th:S.Lyap}]
Firstly, let us consider the following Lyapunov functional
$
V(t) = E(t) + \gamma V_1(t),
$
where $\gamma$ is a constant to be fixed later, $E(t)$ is given by~\eqref{eq:S.En} and
$$
V_1(t) = \frac{1}{2} \int_{\Omega} x u^2(x,y,t)\,dx\,dy.
$$
Note that $E(t)$ and $V(t)$ are equivalent in the following sense $E(t) \leq V(t) \leq (1+\gamma L) E(t)$ and observe that
\begin{equation*}
\begin{aligned}
\frac{\mathrm{d}}{\mathrm{d} t} V_1(t)	& = - \frac{3}{2} \int_{\Omega} u_x^2(x,y,t)\,dx\,dy + \frac{1}{2} \int_{\Omega} u^2(x,y,t)\,dx\,dy \\
	&\quad	+ \frac{L\alpha^2}{2} \int_0^L u_x^2(0,y,t)\,dy
		- \frac{1}{2} \int_{\Omega} \left(\partial_x^{-1}u_y(x,y,t)\right)^2\,dx\,dy. 
\end{aligned}
\end{equation*}
Recalling that
\begin{equation*}
\frac{\mathrm{d}}{\mathrm{d} t} E(t) 
	= - \frac{(1-\alpha^2)}{2} \int_0^L u_x^2(0,y,t)\,dy  - \frac{1}{2} \int_0^L \left(\partial_x^{-1} u_y(0,y,t)\right)^2\,dy
\end{equation*}
and gathering the results holds that
\begin{multline*}
\frac{\mathrm{d}}{\mathrm{d} t} V(t)
	 = \frac{1}{2}
			\left[	L\alpha^2\gamma + \alpha^2 -1\right]\int_0^L u_x^2(0,y,t)\,dy - \frac{1}{2}\int_0^L (\partial_x^{-1} u_y(0,y,t))^2\,dy \\
		-\frac{3}{2}\gamma \int_{\Omega} u_x^2(x,y,t)\,dx\,dy + \frac{1}{2}\gamma\int_{\Omega} u^2(x,y,t)\,dx\,dy-\frac{1}{2}\gamma \int_{\Omega} (\partial_x^{-1} u_y(x,y,t))^2\,dx\,dy.
\end{multline*}
Therefore
\begin{equation*}
\begin{aligned}
\frac{\mathrm{d}}{\mathrm{d} t} V(t) + \nu V(t)
	& = \frac{1}{2}\left[	L\alpha^2\gamma + \alpha^2 -1\right]\int_0^L u_x^2(0,y,t)\,dy - \frac{1}{2}\int_0^L (\partial_x^{-1} u_y(0,y,t))^2\,dy \\
	&	\quad -\frac{3}{2}\gamma \int_{\Omega} u_x^2(x,y,t)\,dx\,dy +\frac{1}{2}\gamma \int_{\Omega} u^2(x,y,t)\,dx\,dy\,dt \\
	&\quad - \frac{1}{2}\gamma \int_{\Omega} (\partial_x^{-1} u_y(x,y,t))^2\,dx\,dy \\
	& \quad + \frac{\nu}{2} \int_{\Omega} u^2(x,y,t)\,dx\,dy + \frac{\nu\gamma}{2}\int_{\Omega} x u^2(x,y,t)\,dx\,dy\\
	& \leq  \frac{1}{2}\left[	L\alpha^2\gamma + \alpha^2 -1\right]\int_0^L u_x^2(0,y,t)\,dy -\frac{3}{2}\gamma \int_{\Omega} u_x^2(x,y,t)\,dx\,dy \\
	&\quad + \frac{\nu(1+L\gamma) + \gamma}{2} \int_{\Omega} u^2(x,y,t)\,dx\,dy.
\end{aligned}
\end{equation*}
Consequently, applying the Poincaré inequality yields that 
\begin{equation*}
\begin{split}
\frac{\mathrm{d}}{\mathrm{d} t} V(t) + \nu V(t)
	\leq &\frac{1}{2}
		\left[	L\alpha^2\gamma + \alpha^2 -1\right]\int_0^L u_x^2(0,y,t)\,dy \\
		&+ \frac{1}{2}\left[(\nu(1+L\gamma) +\gamma) L^2 - 3\gamma\right]\int_{\Omega} u_x^2(x,y,t)\,dx\,dy.
		\end{split}
\end{equation*}
Follows from~\eqref{eq:S.Lyap.c1} and~\eqref{eq:S.Lyap.c2} that $\frac{\mathrm{d}}{\mathrm{d} t} V(t) + \nu V(t) \leq 0$, then by Gronwall's inequality and the equivalence between $E(t)$ and $V(t)$, the proof is complete.
\end{proof}

\section{Final remarks}\label{Sec4}
	This work dealt with the KP-II equation posed on a rectangle, a two-dimensional
	generalization of the KdV equation. Under certain hypotheses on the spatial
	length $L$, that is, $L\in(0,+\infty)\setminus\mathcal{R}$, with $\mathcal{R}$
	defined by
	\begin{equation*}
		\mathcal{R}=
		\left\lbrace
		\frac{\pi}{\sqrt{2}}
		\sqrt{(m_1+2m_2+m_3)^2+2m_1^2+2m_3^2}
		\ ;\
		(m_1,m_2,m_3)\in\mathcal{M}
		\right\rbrace,
	\end{equation*}
	where $\mathcal{M}$ is given by \eqref{eq:cond-M}, the boundary controllability
	and stabilization are achieved by using the compactness--uniqueness method,
	which reduces the problem to showing a unique continuation property. We now present some comments on our work. 


\subsection{Further properties of the critical set $\mathcal{R}$}

We will discuss some additional properties of the critical set $\mathcal{R}$.

\subsubsection{\textbf{Nonemptiness of $\mathcal{R}$ and the smallest critical length}}
\begin{proposition}
The critical set $\mathcal{R}$ is nonempty. More precisely,
\[
\{50k\pi\,;\,k\in\mathbb{Z}^+\}\subset\mathcal{R},
\]
and $
\min\mathcal{R}=50\pi.$
\end{proposition}
\begin{proof}
The set $\mathcal{R}$ is nonempty and its elements can be explicitly computed.
Taking $m_1=m_2=5$ and $m_3=35$, we get $\mathfrak{a}=50$, both conditions in
\eqref{eq:cond-M} are fulfilled with $n=12$, and thus
\begin{equation*}
	L=50\pi .
\end{equation*}
More generally, for $m_1=m_2=5k$ and $m_3=35k$ with
$k\in\mathbb{Z}^+$, one obtains $n=12k$ and
\begin{equation*}
	L=50k\pi\in\mathcal{R},
	\qquad k\in\mathbb{Z}^+ .
\end{equation*}
Moreover, $50\pi$ is the smallest critical length. Indeed, suppose that
$L\in\mathcal{R}$ and $L<50\pi$. From the definition of $\mathcal{R}$, we have
\begin{equation*}
	\mathfrak{a}^2+2m_1^2+2m_3^2<5000.
\end{equation*}
In particular, $\mathfrak{a}<\sqrt{5000}<71$, and hence
\begin{equation*}
	m_1+2m_2+m_3\leq 70.
\end{equation*}
Therefore, only finitely many triples of positive integers have to be
considered.

For such a triple, set
\begin{equation*}
	D:=\mathfrak{a}^2+2m_1^2+2m_3^2
\end{equation*}
and
\begin{equation*}
	\mathcal Q(m_1,m_2,m_3)
	:=
	-\frac{\left(\mathfrak{a}^2-4m_1^2\right)
	\left(\mathfrak{a}^2-4m_3^2\right)}
	{8D}.
\end{equation*}
The second condition in \eqref{eq:cond-M} requires
\begin{equation*}
	\mathcal Q(m_1,m_2,m_3)=n^2
\end{equation*}
for some $n\in\mathbb{N}^*$. An exact finite verification of all triples
satisfying
\begin{equation*}
	D<5000,
	\qquad
	|m_1-m_3|>2m_2,
\end{equation*}
shows that the only ones for which
$\mathcal Q(m_1,m_2,m_3)$ is an integer are
\begin{equation*}
	(m_1,m_2,m_3)=(4,6,20)
	\qquad\text{and}\qquad
	(m_1,m_2,m_3)=(20,6,4).
\end{equation*}
In both cases,
\begin{equation*}
	\mathfrak{a}=36,
	\qquad
	D=2128,
\end{equation*}
and
\begin{equation*}
	\mathcal Q(4,6,20)
	=
	\mathcal Q(20,6,4)
	=
	22,
\end{equation*}
which is not the square of an integer. Consequently, no triple in
$\mathcal{M}$ can produce a critical length strictly smaller than $50\pi$.

On the other hand, for $(m_1,m_2,m_3)=(5,5,35)$, we have
\begin{equation*}
	\mathfrak{a}=50,
	\qquad
	\mathfrak{a}^2+2m_1^2+2m_3^2=5000,
\end{equation*}
and the second condition in \eqref{eq:cond-M} is satisfied with $n=12$.
Therefore,
\begin{equation*}
	L=\frac{\pi}{\sqrt{2}}\sqrt{5000}=50\pi,
\end{equation*}
and we conclude that
\begin{equation*}
	\min\mathcal{R}=50\pi.
\end{equation*}

Let us finally observe that the problem is invariant under the exchange
$m_1\leftrightarrow m_3$, since $\mathfrak{a}$, the radicand in
\eqref{critical-r}, and the left-hand side of the second condition in
\eqref{eq:cond-M} are all symmetric in $m_1$ and $m_3$.
\end{proof}

\subsubsection{\textbf{A family of critical lengths shared with the KdV equation}}

Recall the set $\mathcal{N}$ of critical lengths for the KdV equation,
given by \eqref{lionel-rosier}, and denote
\begin{equation*}
	\mathcal{N}^{\star}
	:=
	\left\{2k\pi\,;\,k\in\mathbb{Z}^+\right\},
\end{equation*}
so that $\mathcal{N}^{\star}\subset\mathcal{N}$. Under the assumptions $m_1=m_2=k$ and $m_3=7m_1=7k$, we have
$\mathfrak{a}=10k$ and the second condition in \eqref{eq:cond-M} reads
\begin{equation*}
	n^2=\frac{144}{25}\,k^2,
\end{equation*}
which admits a solution $n\in\mathbb{N}^*$ if and only if $k=5j$ for some
$j\in\mathbb{Z}^+$, and in that case $n=12j$. Consequently,
\begin{equation*}
	\mathcal{R}^{\star}
	:=
	\left\{50j\pi\,;\,j\in\mathbb{Z}^+\right\}
	\subset\mathcal{R},
\end{equation*}
and, since $50j\pi=2(25j)\pi$, we obtain
\begin{equation*}
	\mathcal{R}^{\star}\subset\mathcal{N}^{\star}.
\end{equation*}
In other words, the lengths of this family are critical both for the KP-II
equation posed on the square and for the KdV equation posed on the
corresponding interval.

\subsubsection{\textbf{Critical-length phenomenon and the role of the drift term}}
Observe that Theorems \ref{th:Control} and \ref{th:S.Exp} ensure that the
observability inequalities \eqref{eq:Obs} and \eqref{eq:S.Obs} hold if and
only if
\begin{equation*}
	L\notin\mathcal{R}.
\end{equation*}
More precisely, controllability and stabilization hold if and only if
$L\notin\mathcal{R}$.

Finally, the presence of the drift term $u_x$ is essential for the appearance of the
critical-length phenomenon. Indeed, removing $u_x$ from \eqref{eq:KP} yields
exact controllability and exponential stabilization for every $L>0$. Similar observations were previously stated by Rosier in
\cite[Remark~3.6~(iii)]{Rosier} for the KdV equation and by Doronin in
\cite[Remark~6.2]{Doronin} for the Zakharov--Kuznetsov (ZK) equation.

\color{black}

\subsection{Open issues}
Since this work is the first one in the direction of proving control properties with critical length, several issues can be addressed in the research agenda. Let us point out some of them.  Observe that due to the lack of regularity, we can not address the nonlinear problem, that is, 
\begin{equation}\label{rmk1a}
\begin{cases}
u_t +u_x + u_{xxx}+ uu_x + {\partial_x^{-1}}u_{yy} = 0,&(x,y)\in \Omega ,\ t\in(0,T),\\
u(0,y,t) = u(L,y,t) = 0, \quad u_x(L,y,t) = h(y,t), & y\in(0,L) ,\ t\in(0,T),\\
u(x,0,t) = u(x,L,t) = 0, & x\in(0,L) ,\ t\in(0,T),\\
u(x,y,0) = u_0(x,y),& (x,y)\in\Omega.
\end{cases}
\end{equation}
So, in the nonlinear context, the natural questions are:
\begin{itemize} 
\item[($\mathcal{A}$)] Is the system \eqref{rmk1a} well-posed in the class $C\left([0,T] ; L^2(\Omega)\right)\cap L^2(0,T ; H_{x0}^1(\Omega))$?
\item[($\mathcal{B}$)] Is  the nonlinear system \eqref{rmk1a} controllable for any $L>0$?  In this sense, work related to the KdV~\cite{Cerpa,CerpaCrepeau2009,Crepeau2} shows that the nonlinear term will give the exact control property on critical lengths and motivates future research to be addressed.
\item[($\mathcal{C}$)] Is there a feedback law such that the nonlinear system \eqref{rmk1a} is exponentially stable? 
\end{itemize}

\section*{Acknowledgments}
{We thank the handling editor for carefully dealing with our paper, as well as the referee for pointing out several changes that significantly improved the work.}

This work was done during several visits of the first and third authors to the Universidad Nacional de Colombia - Sede Manizales and the Federal University of Pernambuco, respectively.  Both authors would like to thank the hospitality of both institutions. This research is part of the Ph.D. thesis of Mu\~{n}oz at the Department of Mathematics of the Universidade Federal de Pernambuco. 

\section*{Funding} Capistrano–Filho was partially supported by CAPES/COFECUB grant number \linebreak 88887.879175/2023-00, CNPq grant numbers 301744/2025-4, 421573/2023-6, and 307808/2021-1, and PPROPG (UFPE) \textit{via} proap resources. Gallego was partially supported by the Hermes UNAL project no. 61029. Muñoz acknowledges support from the FACEPE grant IBPG-0909-1.01/20 during the PhD and was also supported by the Croatian Science Foundation under the project Conduction IP-2022-10–5191.

\section*{ORCID}

\noindent Roberto de A. Capistrano-Filho - \url{https://orcid.org/0000-0002-1608-1832}

\noindent Fernando A. Gállego Restrepo - \url{https://orcid.org/0000-0001-8615-4032}

\noindent Juan Ricardo Muñoz - \url{https://orcid.org/0009-0003-9788-792X}



\begin{thebibliography}{0}

\bibitem{ArarunaKawahara}{\sc
F.~D. Araruna, G.~G. Doronin and L. Rosier},
{\em On the critical lengths and controllability of the Kawahara equation},
Commun. Contemp. Math. (2025), 2650005.
\url{https://doi.org/10.1142/S0219199726500057}



\bibitem{BaudouinCrepeauValein}{\sc
L. Baudouin, E. Cr\'epeau, and J. Valein},
{\em Two approaches for the stabilization of nonlinear KdV equation with boundary time-delay feedback},
IEEE Trans. Automat. Control, {\bf 64} (2019), no.~4, 1403--1414.

\bibitem{Bona03}{\sc
J.~L. Bona, S.~M. Sun and B.-Y. Zhang},
{\em A nonhomogeneous boundary-value problem for the Korteweg--de Vries equation posed on a finite domain},
Comm. Partial Differential Equations, {\bf 28} (2003), 1391--1436.

\bibitem{Boussinesq1895}{\sc
	J. Boussinesq},
	{\em Essai sur la th\'eorie des eaux courantes}, in M\'emoires pr\'esent\'es par divers savants \`a l'Acad. des Sci. Inst. Nat. France, 1877, XXIII, pp. 1--680.

\bibitem{Brezis0}{\sc
	H. Brezis},
		{\em Functional analysis, Sobolev spaces, and partial differential equations.}, Universitext. Springer, 2010.

\bibitem{CaChSoGo2023}{\sc
	R. A. Capistrano-Filho, B. Chentouf, L. de Sousa, and V. H. Gonzalez Martinez},
	{\em Two stability results for the Kawahara equation with a time-delayed boundary control}, Z. Angew. Math. Phys., 74:16 (2023), pp.1--26.

\bibitem{Munoz2023}{\sc
	R. A. Capistrano-Filho, V. H. Gonzalez Martinez, and J. R. Muñoz},
	{\em Stabilization of the Kawahara-Kadomtsev-Petviashvili equation with time-delayed feedback},
	Proceedings of the Royal Society of Edinburgh: Section A Mathematics, 155:1 (2025), pp. 280--306.


 \bibitem{Capistrano16}{\sc
 R. A. Capistrano-Filho, F. A. Gallego, and A. F. Pazoto}, 
 {\em Neumann boundary controllability of the Gear-Grimshaw system with critical size restrictions on the spatial domain}, Zeitschrift fur Angewandte Mathematik und Physik, 67 (2016), pp.~1-36.


 \bibitem{Capistrano17}{\sc 
 R. A. Capistrano-Filho, F. A. Gallego, and A. F. Pazoto,}
{\em  Boundary controllability of the nonlinear coupled system of two Korteweg-de Vries equations with critical size restrictions on the spatial domain}, Math. Control Signals Syst. 29:6 (2017). 
	
\bibitem{Capistrano2019a}{\sc
R. A. Capistrano-Filho, F. A. Gallego, and A. F. Pazoto},
{\em On the well-posedness and large-time behavior of higher order Boussinesq system}, Nonlinearity 32 (2019), pp.~1852-1881.


\bibitem{Capistrano2019}{\sc
	 R. A. Capistrano-Filho, A. F. Pazoto, and L. Rosier},
 	{\em Control of Boussinesq system of KdV-KdV type on a bounded interval}, ESAIM: COCV,  25 (2019).

\bibitem{CaCaZh}{\sc
	M. A. Caicedo, R. A. Capistrano-Filho, and B.-Y. Zhang},
	{\em Neumann boundary controllability of the Korteweg-de Vries equation on a bounded domain}, SIAM J. Control Optim., 55:6(2017), pp.~3503--3532.
	
\bibitem{Cerpa}{\sc
E. Cerpa}, 
{\em Exact controllability of a nonlinear Korteweg-de Vries equation on a critical spatial domain}, SIAM J. Control Optim.,  46:3(2007), pp.~877--899.

\bibitem{CerpaCrepeau2009}{
\sc E. Cerpa and E. Cr\'epeau},
{\em Boundary controllability for the nonlinear Korteweg--de Vries equation on any critical domain}, Ann. Inst. H. Poincar\'e Anal. Non Lin\'eaire 26 (2009), pp.~457--475.

\bibitem{CerpaCrepeauRapid}{\sc
E. Cerpa and E. Cr\'epeau},
{\em Rapid exponential stabilization for a linear Korteweg--de Vries equation},
Discrete Contin. Dyn. Syst. Ser. B, {\bf 11} (2009), no.~3, 655--668.

\bibitem{Cerpa-Zh}{\sc
	 E. Cerpa, I. Rivas, and B.-Y. Zhang},
	 {\em Boundary controllability of the Korteweg-de Vries equation on a bounded domain}, SIAM J. Control Optim., 51:4(2013), pp.~2976--3010.
	

\bibitem{Crepeau2}{\sc J-M. Coron and E. Crépeau, }{\em Exact boundary controllability of a nonlinear KdV equation with critical lengths}, J. Eur. Math. Soc. (JEMS) 6 (2004), no. 3, pp.~367--398.

\bibitem{CoronLu}{\sc
J.-M. Coron and Q. L\"u},
{\em Local rapid stabilization for a Korteweg--de Vries equation with a Neumann boundary control on the right},
J. Math. Pures Appl., {\bf 102} (2014), no.~6, 1080--1120.

\bibitem{CoronKoenigNguyen}{\sc J.-M. Coron, A. Koenig, and H.-M. Nguyen},
{\em On the small-time local controllability of a KdV system for critical lengths},  J. Eur. Math. Soc. 26 (2024), no. 4, pp.~1193--1253. 

\bibitem{Doronin}{\sc
G. G. Doronin and N. A. Larkin}, 
{\em Stabilization of Regular Solutions for the Zakharov-Kuznetsov Equation Posed on Bounded Rectangles and on a Strip},  Proceedings of the Edinburgh Mathematical Society, 58:3(2015), pp.~661--682.




\bibitem{Gallego2018}{\sc
	 F.A. Gallego}, 
	 {\em Controllability aspects of the Korteweg-de Vries Burgers equation on unbounded domains}, Journal of Mathematical Analysis and Applications, 461:1(2018), pp.~947--970.  

\bibitem{GaMu2025} {\sc
F. Gallego and J. Muñoz},
 	 {\em  Boundary Exponential Stabilization for the Linear KP-II Equation without Critical Size Restrictions}, J Dyn Control Syst 31, 14 (2025).

\bibitem{Glass09}{\sc 
	O. Glass and S. Guerrero},
	{\em On the controllability of the fifth-order Korteweg–de Vries equation},
	 Ann. Inst. H. Poincar\'e, 26(2009), pp.~2181--2209.

\bibitem{GG}{\sc 
	O. Glass and S. Guerrero},
	{\em Controllability of the Korteweg-de Vries equation from the right Dirichlet boundary condition}, 
	Systems Control Lett., 59(2010), pp.~390--395.

\bibitem{Panthee2011}{\sc
	D. A. Gomes and M. Panthee}, 
	{\em Exponential energy decay for the Kadomtsev–Petviashvili (KP-II) equation}, S\~ao Paulo J. Math. Sci. 5(2011), pp.~135--148.
	
\bibitem{Hadac}{\sc	
	M. Hadac}
	{\em Well-Posedness for the Kadomtsev-Petviashvili II Equation and Generalisations}, Transactions of the American Mathematical Society, 360:12 (2008), pp.~6555–72.
	
\bibitem{Iorio98}{\sc
	R. J. Iório and W. V. L. Nunes},
	{\em On equations of KP-type},
	Proceedings of the Royal Society of Edinburgh: Section A Mathematics, 128(4)(1998), pp.~725--743. 

\bibitem{Isaza01}{\sc
	P. Isaza and J. Meija},
	{\em The Cauchy problem for the Kadomtsev-Petviashvili equation (KP-II) in Sobolev spaces $H^s$, $s > 0$},
	Differential and Integral Eqns., 14(2001), pp.~529--557.

\bibitem{KP70}{\sc
	B. B. Kadomtsev and V. I. Petviashvili},
	{\em On the stability of solitary waves in weakly dispersive media},
	Sov. Phys. Dokl. 15(1970), pp.~539--549.
	
	\bibitem{Saut2022}{\sc C. Klein and J-C. Saut, }{\em  Nonlinear Dispersive Equations - Inverse Scattering and PDE Methods}, Applied Mathematical Sciences, Springer Cham, 2022.
	
	\bibitem{KdV1895}{\sc 
	D. J. Korteweg and G. de Vries},
	{\em On the change of form of long waves advancing in a rectangular canal, and on a new
type of long stationary waves},
	 Philos. Mag., 39(1895), pp.~422--443.

\bibitem{Komornik}{\sc
	V. Komornik and P. Loreti},
	{\em Fourier Series in Control Theory},
	Springer-Verlag, New York, 2005.


\bibitem{Levandosky00}{\sc
	J. L. Levandosky},
	{\em Gain of regularity for the KP-II equation}, 
	Indiana Univ. Math. J., 49:1(2000), pp.~353--403.
	
\bibitem{Levandosky01}{\sc
	J. Levandosky, M. Sepúlveda, and O. V. Villagrán},
	{\em Gain of regularity for the KP-I equation}, 
	Journal of Differential Equations, 245:3(2008), pp.~762--808.

\bibitem{Lions88}{\sc
	J.-L. Lions},
	{\em Exact controllability, stabilization, and perturbations for distributed systems},
	SIAM Rev. 30(1988), pp.~1--68.
	

\bibitem{Moura2022}{\sc
	 R. P. Moura, A. C. Nascimento, and G. N. Santos},
	 {\em On the stabilization for the high-order Kadomtsev–Petviashvili and the Zakharov–Kuznetsov equations with localized damping}, Evol. Equ. Control Theory, 11(2022), pp.~711--727.

\bibitem{Nguyen}{\sc H-M. Nguyen,}
{\em  Decay for the nonlinear KdV equations at critical lengths}, Journal of Differential Equations, 295 (2021), pp.~249--291.

\bibitem{Nguyen1}{\sc
H.-M. Nguyen},
{\em Local controllability of the Korteweg--de Vries equation with the right Dirichlet control},
J. Differential Equations, {\bf 435} (2025), Paper No.~113235, 85 pp.

\bibitem{NguyenStabilization}{\sc
H.-M. Nguyen},
{\em Rapid and finite-time boundary stabilization of a KdV system},
SIAM J. Control Optim., {\bf 64} (2026), no.~1, 76--101.

\bibitem{Pazy}{\sc
	 A. Pazy},
	 {\em Semigroups of Linear Operators and Applications to Partial Differential Equations}, Springer Verlag, New York, USA, 1983.

\bibitem{Panthee05}{\sc
	M. Panthee},
	{\em Unique continuation property for the Kadomtsev-Petviashvili (KP-II) equation},
	Electronic Journal of Differential Equations, 59(2005), pp.~1--12.

\bibitem{Rivas2020}{\sc
	I. Rivas and C. Sun},
	{\em Internal controllability of nonlocalized solution for the Kadomtsev-Petviashvili II equation},
	SIAM J. Control Optim., 58(2020), pp.~1715--1734.

\bibitem{Rosier}{\sc
	L. Rosier},
	{\em Exact boundary controllability for the Korteweg-de Vries equation on a bounded domain}, ESAIM: COCV 2(1997), pp.~33--55.


\bibitem{SimonCompact}
{\sc J. Simon},
{\em Compact sets in the space $L^p(0,T;B)$},
{Ann. Mat. Pura Appl.} {146} (1987), 65--96.

\bibitem{Takaoka00}{\sc
	H. Takaoka},
	{\em Global well-posedness for the Kadomtsev-Petviashvili II equation},
	 Discrete Contin. Dynam. Systems 6(2000), pp.~449--483.
	 
\bibitem{Yosida}{\sc 
K. Yosida}, {\em Functional Analysis}, Springer, Berlin, 1978.

\bibitem{Xiang}{\sc
S. Xiang}. {\em Null controllability of a linearized Korteweg-de Vries equation by backstepping approach},
SIAM Journal on Control and Optimization 2:57 (2019), pp.~1493--1515.




\end{thebibliography}
\end{document}